\documentclass[a4paper,12pt,twoside]{article}
\usepackage{graphicx}\usepackage{a4wide}
\usepackage{hyperref}
\usepackage{amssymb,amsmath,amsthm}
\usepackage{setspace}
\usepackage{tikz}
\usepackage{psfrag}
\usetikzlibrary{shadows}
\usetikzlibrary{shapes,arrows}
\usepackage{color, colortbl}
\usepackage{arydshln}
\usepackage{float}
\usepackage{ulem}
\usepackage{subcaption}
\usepackage[font=normalsize,labelfont=bf]{caption}
\usepackage[font=normalsize,labelfont=bf]{subcaption}
\usepackage{cleveref}
\usepackage{upgreek,xspace}
\usepackage{commath}

\makeatletter
\newcommand*\wt[2][0.2ex]{%
	\begingroup
	\mathchoice{\wt@helper{#1}{#2}{\displaystyle}{\textfont}}
	{\wt@helper{#1}{#2}{\textstyle}{\textfont}}
	{\wt@helper{#1}{#2}{\scriptstyle}{\scriptfont}}
	{\wt@helper{#1}{#2}{\scriptscriptstyle}{\scriptscriptfont}}%
	\endgroup
	#2%
}
\newcommand*\wt@helper[4]{%
	\def\currentfont{\the#41}%
	\def\currentskewchar{\char\the\skewchar\currentfont}%
	\setbox\tw@\hbox{\currentfont#2\currentskewchar}%
	\dimen@ii\wd\tw@
	\setbox\tw@\hbox{\currentfont#2{}\currentskewchar}%
	\advance\dimen@ii-\wd\tw@
	\rlap{\raisebox{-#1}{$\m@th#3\kern\dimen@ii\widetilde{\phantom{#2}}$}}%
}
\makeatother

\newcommand{\bm}[1]{\text{\boldmath $#1$\unboldmath}}

\newcommand{\Rd}{\mathbb{R}^d}
\newcommand{\grad}{\bm{\nabla}}
\newcommand{\bsigma}{\bm{\sigma}}
\newcommand{\bu}{\bm{u}}
\newcommand{\bb}{\bm{b}}
\newcommand{\bt}{\bm{t}}
\newcommand{\bn}{\bm{n}}
\newcommand{\bE}{\bm{E}}
\newcommand{\beps}{\bm{\varepsilon}}
\newcommand{\bD}{\mathbf{D}}

\newcommand{\acc}{\ddot{\bu}}

\newcommand{\bv}{\bm{v}}
\newcommand{\bw}{\bm{w}}
\newcommand{\Hone}{\mathcal{H}^1}
\newcommand{\Hd}{\mathcal{H}^1_D}
\newcommand{\Hz}{\mathcal{H}^1_0}

\newcommand{\bx}{\bm{x}}
\newcommand{\bxi}{\bm{\xi}}

\newcommand{\nel}{n_{el}}
\newcommand{\nen}{n_{en}}
\newcommand{\nN}{n_{n}}
\newcommand{\ndof}{n_{d}}
\newcommand{\bN}{\bm{N}}

\newcommand{\refElem}{\widehat{\Omega}}
\newcommand{\isomap}{\bm{\varphi}^e}
\newcommand{\bJ}{\mathbf{J}}
\newcommand{\isoJ}{\bJ_\isomap}
\newcommand{\isoJF}{\bJ_\Gamma^e}

\newcommand{\bB}{\mathbf{B}}
\newcommand{\bM}{\mathbf{M}}
\newcommand{\bK}{\mathbf{K}}
\newcommand{\bU}{\mathbf{U}}
\newcommand{\bF}{\mathbf{F}}

\newcommand{\bI}{\mathbf{I}}
\newcommand{\nr}{n_{r}}

\newcommand{\bKll}{\bK_{ll}}
\newcommand{\bKls}{\bK_{ls}}
\newcommand{\bKsl}{\bK_{sl}}
\newcommand{\bKss}{\bK_{ss}}
\newcommand{\bMll}{\bM_{ll}}
\newcommand{\bMls}{\bM_{ls}}
\newcommand{\bMsl}{\bM_{sl}}
\newcommand{\bMss}{\bM_{ss}}
\newcommand{\balpha}{\boldsymbol{\alpha}}
\newcommand{\bPhi}{\mathbf{\Phi}}

\newcommand{\bmu}{\bm{\mu}}

\newcommand{\mM}{\mathcal{M}}
\newcommand{\bmM}{\bm{\mM}}

\newcommand{\np}{n_p}

\renewcommand{\bf}{\bm{f}}

\newcommand{\fb}{\mathbf{b}}

\newcommand{\fx}{\mathbf{x}}

\newcommand{\fA}{\mathbf{A}}

\newcommand{\fF}{\mathbf{F}}

\newcommand{\fK}{\mathbf{K}}

\newcommand{\fM}{\mathbf{M}}

\newcommand{\fU}{\mathbf{U}}

\newcommand{\T}[1]{{#1}^{\sf  T}}

\DeclareMathAlphabet\mathbfcal{OMS}{cmsy}{b}{n}

\newenvironment{keywords}{\begin{quote}\emph{\textbf{Keywords:}}}{\end{quote}}

\theoremstyle{definition}

\newtheorem{remark}{Remark}

\makeatletter
\pgfkeys{/pgf/.cd,
	parallelepiped offset x/.initial=2mm,
	parallelepiped offset y/.initial=2mm
}
\pgfdeclareshape{parallelepiped}
{
	\inheritsavedanchors[from=rectangle] 
	\inheritanchorborder[from=rectangle]
	\inheritanchor[from=rectangle]{north}
	\inheritanchor[from=rectangle]{north west}
	\inheritanchor[from=rectangle]{north east}
	\inheritanchor[from=rectangle]{center}
	\inheritanchor[from=rectangle]{west}
	\inheritanchor[from=rectangle]{east}
	\inheritanchor[from=rectangle]{mid}
	\inheritanchor[from=rectangle]{mid west}
	\inheritanchor[from=rectangle]{mid east}
	\inheritanchor[from=rectangle]{base}
	\inheritanchor[from=rectangle]{base west}
	\inheritanchor[from=rectangle]{base east}
	\inheritanchor[from=rectangle]{south}
	\inheritanchor[from=rectangle]{south west}
	\inheritanchor[from=rectangle]{south east}
	\backgroundpath{
		\southwest \pgf@xa=\pgf@x \pgf@ya=\pgf@y
		\northeast \pgf@xb=\pgf@x \pgf@yb=\pgf@y
		\pgfmathsetlength\pgfutil@tempdima{\pgfkeysvalueof{/pgf/parallelepiped offset x}}
		\pgfmathsetlength\pgfutil@tempdimb{\pgfkeysvalueof{/pgf/parallelepiped offset y}}
		\def\ppd@offset{\pgfpoint{\pgfutil@tempdima}{\pgfutil@tempdimb}}
		\pgfpathmoveto{\pgfqpoint{\pgf@xa}{\pgf@ya}}
		\pgfpathlineto{\pgfqpoint{\pgf@xb}{\pgf@ya}}
		\pgfpathlineto{\pgfqpoint{\pgf@xb}{\pgf@yb}}
		\pgfpathlineto{\pgfqpoint{\pgf@xa}{\pgf@yb}}
		\pgfpathclose
		\pgfpathmoveto{\pgfqpoint{\pgf@xb}{\pgf@ya}}
		\pgfpathlineto{\pgfpointadd{\pgfpoint{\pgf@xb}{\pgf@ya}}{\ppd@offset}}
		\pgfpathlineto{\pgfpointadd{\pgfpoint{\pgf@xb}{\pgf@yb}}{\ppd@offset}}
		\pgfpathlineto{\pgfpointadd{\pgfpoint{\pgf@xa}{\pgf@yb}}{\ppd@offset}}
		\pgfpathlineto{\pgfqpoint{\pgf@xa}{\pgf@yb}}
		\pgfpathmoveto{\pgfqpoint{\pgf@xb}{\pgf@yb}}
		\pgfpathlineto{\pgfpointadd{\pgfpoint{\pgf@xb}{\pgf@yb}}{\ppd@offset}}
	}
}
\makeatother
\graphicspath{{Figures/}}

\begin{document}
\title{Nonintrusive reduced order model for parametric solutions of inertia relief problems}

\author{
\renewcommand{\thefootnote}{\arabic{footnote}}
			  F. Cavaliere\footnotemark[1]\textsuperscript{ \ ,}\footnotemark[2]\textsuperscript{ \ ,}\footnotemark[3] \ ,
			  S. Zlotnik\footnotemark[1]\textsuperscript{ \ ,}\footnotemark[2], \
			  R. Sevilla\footnotemark[3], \
             X.Larr\'{a}yoz\footnotemark[4] \ and
             P. D\'{i}ez\footnotemark[1]\textsuperscript{ \ ,}\footnotemark[2]\
}

\date{\today}
\maketitle

\renewcommand{\thefootnote}{\arabic{footnote}}

\footnotetext[1]{Laboratori de C\`alcul Num\`eric (LaC\`aN), ETS de Ingenieros de Caminos, Canales y Puertos, Universitat Polit\`ecnica de Catalunya, Barcelona, Spain}
\footnotetext[2]{International Centre for Numerical Methods in Engineering, CIMNE, Barcelona, Spain}
\footnotetext[3]{Zienkiewicz Centre for Computational Engineering, College of Engineering, Swansea University, Wales, UK}
\footnotetext[4]{Centro T\'{e}cnico de SEAT S.A, Autov\'{i}a A-2, km 85, Martorell, Spain
\vspace{5pt}\\
Corresponding author: Fabiola Cavaliere. \textit{E-mail:} \texttt{fabiola.cavalere1@upc.edu}
}

\begin{abstract}
	The Inertia Relief (IR) technique is widely used by industry and produces equilibrated loads allowing to analyze unconstrained systems without resorting to the more expensive full dynamic analysis. The main goal of this work is to develop a computational framework for the solution of unconstrained parametric structural problems with IR and the Proper Generalized Decomposition (PGD) method. First, the IR method is formulated in a parametric setting for both material and geometric parameters. A reduced order model using the {\textit{encapsulated PGD}} suite is then developed to solve the parametric IR problem, circumventing the so-called {\textit{curse of dimensionality}}.
		With just one offline computation, the proposed PGD-IR scheme provides a {\textit{computational vademecum}} that contains all the possible solutions for a pre-defined range of the parameters. The proposed approach is nonintrusive and it is therefore possible to be integrated with commercial FE packages. 
		The applicability and potential of the developed technique is shown using a three dimensional test case and a more complex industrial test case. The first example is used to highlight the numerical properties of the scheme, whereas the second example demonstrates the potential in a more complex setting and it shows the possibility to integrate the proposed framework within a commercial FE package. In addition, the last example shows the possibility to use the generalized solution in a multi-objective optimization setting.
\end{abstract}
\begin{keywords}
Proper generalized decomposition, inertia relief, nonintrusive, reduced order model, shape optimization
\end{keywords}
\section{Introduction}\label{secIntro}
Unconstrained structures are widespread in the automotive, aerospace and naval industry. As is well known, due to the singularity of the stiffness matrix, conventional static analyses cannot be performed if the system undergoes rigid body motions. At the same time, imposing dummy constraints in order to make a free-body system statically determinate leads to unrealistic reaction forces and, as a consequence, an unrealistic distribution of the internal stresses. The inertia relief (IR) method represents an attractive alternative for solving unconstrained structural problems without resorting to the more expensive full dynamic analysis~\cite{wijker2004mechanical}. The main idea is to counteract the unbalanced applied loads by a set of rigid body accelerations, the latter providing body forces which are distributed over the structure in such a way that the applied forces are equilibrated and the static analysis can be performed. The technique is available into most of the commercial finite element packages and it has been widely used by the industry in different fields~\cite{bisplinghoff1955aeroelasticity,  nelson1977use, kuo1995vehicle, barnett1995closed, baskar1998door, anvari1999automotive,pagaldipti2004influence, liao2011study, pengqiu2017modified}.

The static global stiffness analysis of a body in white (BIW) is a common example that involves the computational simulation of an unconstrained configuration using the IR method. The BIW global stiffness plays a significant role in the design process of a car. An important challenge of this problem is the number of parameters (e.g. geometry, material) to be considered during the analysis of a BIW. As any change in the material or geometrical characteristics of the car components might have considerable effects on the global behaviour of the structure, the number of simulations that are required to account for the whole range of the involved parameters becomes prohibitive when classical numerical approaches are employed. As a consequence, the possibility to perform parametric studies, shape optimization or inverse identification in the context of a BIW remains a challenge.

A way to circumvent this issue and to reduce the computational complexity of parametric problems is to employ a reduced order model (ROM). In the last two decades, researchers from the most diverse areas of science and engineering have developed different ROM techniques, with the common goal of finding low-order models, described by a reduced order basis, which are able to capture the essential behaviour of a complex system. Well known computational approaches based on this idea are Krylov-based methods \cite{freund2003model}, the reduced basis method~\cite{rozza2007reduced} and the proper orthogonal decomposition (POD)~\cite{berkooz1993proper, feeny1998physical, rowley2005model}. These techniques, known as \textit{a posteriori} methods, first solve the full-order problem for a suitably chosen set of parameters, providing a set of snapshots of the solution. This step, usually referred to as the offline stage, is used to extract the most relevant characteristics of the solution. Then, during the online phase, the solution for any new parameter value can be expressed as a linear combination of the previously computed basis functions. 

A valuable alternative is represented by the proper generalized decomposition (PGD) method~\cite{ammar2006new,chinesta2013pgd, chinesta2013proper, chinesta2014pgd}, which is an \textit{a priori} approach. The main idea behind the PGD method is to consider the parameters as extra coordinates of the problem, increasing the dimensionality of the problem at hand, and to assume that the solution of the high-dimensional problem can be approximated by a separable function. During the offline stage, usually performed by employing high performance computing resources, the PGD algorithm computes on-the-fly a set of basis functions, usually called modes. The PGD approximation depends explicitly on the parameters, so it represents a \textit{computational vademecum} containing the solution for every possible combination of the parameters. In the online stage, the solution can be particularized for any set of the parameters in real-time. The method has been tested in the most diverse fields, such as flow problems~\cite{dumon2010proper,dumon2011proper,leblond2014priori,diez2017generalized,ibanez2017simulating},  thermal problems~\cite{ghnatios2012proper,aguado2015real,huerta2018proper}, solid mechanics~\cite{de2013basis,reiserror}, fracture mechanics~\cite{giner2013proper,garikapati2020proper}, elastic metamaterials and coupled magneto-mechanical problems~\cite{sibileau2018explicit,barroso2020staggered}. 

Despite the wide range of problems where the PGD has shown its potential, the application to geometrically parametrized problems remains particularly challenging, due to the difficulty to obtain a separable expression of the discrete problem. Previous works that dealt with parametric shapes are usually limited to simple geometric dependence
~\cite{chinesta2013pgd, leygue2010first, bognet2012advanced, heuze2016parametric, courard2016integration}. Other authors proposed a technique based on the idea that a parent domain can be associated to the parametric domain in order to introduce the parametric dependency on the geometry in the governing equations~\cite{leygue2010first, zlotnik2015proper}. More recently, another approach was proposed~\cite{sevilla2020solution} in which the control points characterising the NURBS curves or surfaces used in CAD representation are defined as the geometric parameters of the problem.  

One of the main drawbacks of the original PGD approach, when compared to other \textit{a posteriori} approaches, is the intrusivity of its implementation, which precludes its wider application in an industrial context, where commercial software are usually employed. This limitation has motivated the development of nonintrusive implementations of the PGD rationale for solid~\cite{zou2018nonintrusive} and fluid~\cite{tsiolakis2020nonintrusive} mechanics problems. Also motivated by the goal of achieving the nonintrusive applications of the PGD, the \textit{Encapsulated PGD Toolbox} has been recently developed by D\'{i}ez et al.~\cite{diez2018algebraic,diez2019encapsulated}. The toolbox consists of a collection of PGD-based routines that are able to perform algebraic operations for multidimensional tensors. One relevant advantage of the toolbox is that each routine is \textit{encapsulated} and can be used as a black box, enabling the nonintrusive coupling with commercial FE packages. This feature is of major importance for the application of the PGD method in an industrial setting.

This work proposes a nonintrusive PGD-IR method for the solution of an unconstrained structure characterized by material and/or geometric parameters. The proposed parametric IR approach involves a "cascade" application of the PGD method in order to solve three sequential parametric problems, where the parametric solution of one problem is taken as the input of the next parametric problem.  In order to automate the process, an \textit{ad-hoc} solver was implemented that makes use of the \textit{Encapsulated PGD Toolbox}. In the present work, a nonintrusive interaction between the external commercial finite element (FE) software MSC-Nastran and an in-house code implemented in Matlab is considered.

The structure of the remainder of the paper is as follows. Sec.~\ref{secProblStat} presents the problem statement in terms of an elastodynamic boundary value problem. Sec.~\ref{secIR} briefly reviews the idea behind the IR method for a non-parametric problem. The proposed PGD-IR approach is presented in~\ref{secParamIR}, where the "cascade" application of the PGD approach is proposed to solve the PGD-IR problem. Also, the algebraic approach to deal with geometric parameters is detailed. In Sec.~\ref{secNumEx} two numerical examples are used to show the potential of the proposed method. In the first example, a simple linear elastic 3D structure with one material and one geometric parameter are considered to underline the main properties of the developed ROM. The second example uses a more realistic case of a dummy car to demonstrate the nonintrusive interaction with the commercial FE software MSC-Nastran. A multiobjective optimization study is shown which proves that the method can be employed as a fast and reliable tool to guide the designers in the intricate decision-making procedure. Finally, Sec.~\ref{secConcl} summarizes the conclusions of the work that has been presented.

\section{Finite elements formulation of elastodynamic problems}\label{secProblStat}
%
\subsection{Problem statement} \label{secStrongForm}

Let us consider an open bounded domain $\Omega \subset \Rd$, where $d$ is the number of spatial dimensions. The boundary of the domain is assumed to be partitioned into the the disjoint parts $\Gamma_D$ and $\Gamma_N$, where Dirichlet and Neumann boundary conditions are prescribed respectively. The strong form of the elastodynamic problem using the classical Voigt notation~\cite{jacob2007first} can be written as
\begin{equation}
\begin{cases}\label{eqStrongForm}
\rho \acc -  \grad_S^T \bsigma  =  \bb \quad &\text{in} \quad \Omega \times (0,T] \\
\bu = \bu_D  \quad & \text{on} \quad \Gamma_D \times (0,T] \\
\bE^T \bsigma = \bt_N  \quad &\text{on} \quad \Gamma_N \times (0,T], \\
\bu = \bu_0  \quad & \text{in} \quad \Omega \times \{0\} \\
\dot{\bu} = \bv_0  \quad & \text{in} \quad \Omega \times \{0\}
\end{cases}
\end{equation}
where $\bu$ is the displacement field, $\acc$ denotes the acceleration, $\bsigma$ is a vector containing the extensional and shear stress components of the Cauchy stress tensor, $\bb$ is a external body force vector, $T$ is the final time of interest, $\bE$ is a matrix accounting for the normal direction to the boundary, $\bu_D$ and $\bt_N$ are the prescribed displacement and traction vectors on the Dirichlet and Neumann boundaries respectively and $\bu_0$ and $\bv_0$ are the initial position and velocity respectively. In three dimensions, the matrix operator $\grad_S$ and the matrix $\bE$ are given by
\begin{equation} \label{eqVoigtDefs}
\grad_S :=
\begin{bmatrix}
\partial/\partial x_1 & 0 & 0 & \partial/\partial x_2 & \partial/\partial x_3 & 0 \\
0 & \partial/\partial x_2 & 0 & \partial/\partial x_1 & 0 & \partial/\partial x_3 \\
0 & 0 & \partial/\partial x_3 & 0 & \partial/\partial x_1 & \partial/\partial x_2
\end{bmatrix}^T, 
\qquad
\bE := 
\begin{bmatrix}
n_1 & 0 & 0 & n_2 & n_3 & 0\\
0 & n_2 & 0 & n_1 & 0 & n_3 \\
0 & 0 & n_3 & 0 & n_1 & n_2
\end{bmatrix}^T,
\end{equation}
with $\bn$ being the outward unit normal vector to $\partial \Omega$.
For a linear elastic material, the generalized Hooke's law expresses a linear relation between the stress vector, $\bsigma$, and the strain vector, $\beps$, namely
\begin{equation}\label{eqConstLaw}
\bsigma = \bD \beps,
\end{equation}
where $\beps := \grad_S \bu$ and $\bD$ is a symmetric positive definite matrix depending upon the Young modulus, $E$, and the Poisson ratio, $\nu$. In three dimensions
\begin{equation} \label{eqHooke}
\bD := \frac{E}{(1 + \nu)(1 - 2\nu)}
\begin{bmatrix}
1-\nu & \nu & \nu & \\
\nu & 1-\nu & \nu & \mathbf{0}_d \\
\nu & \nu & 1-\nu & \\
& \mathbf{0}_d & & (1-2\nu)/2 \mathbf{I}_d
\end{bmatrix}.
\end{equation}
The weak formulation of the strong form of Eq.~\eqref{eqStrongForm} reads as follows: given $\bu_D$ on $\Gamma_D$ and $\bt_N$ on $\Gamma_N$, find $\bu \in \Hd(\Omega):= \{ \bw \in \Hone(\Omega) \; | \; \bw = \bu_D  \text{ on } \Gamma_D \}$ such that
\begin{equation}\label{eqWeak}
\int_{\Omega} \rho \bv \cdot \ddot{\bu}  d \Omega + \int_{\Omega} \grad_S \bv \cdot \left( \bD \grad_S \bu \right)  d \Omega = \int_{\Omega} \bv \cdot \bb d\Omega  +  \int_{\Gamma_N} \bv \cdot \bt_N d\Gamma,
\end{equation}
for all $\bv \in \Hz(\Omega):= \{ \bw \in \Hone(\Omega) \; | \; \bw = \bm{0} \text{ on } \Gamma_D \}$.
\subsection{Spatial discretization}\label{secDiscr}
A partition of the domain $\Omega$ in a set of $\nel$ disjoint elements $\Omega_e$ is considered. Following the classical isoparametric framework, the approximation of the displacement field is defined in a reference element, $\refElem$, with local coordinates $\bxi$, as
\begin{equation} \label{eqApprox}
\bu(\bxi) \simeq \bu^h(\bxi) := \sum_{j=1}^{\nen} \bU_j N_j(\bxi),
\end{equation}
where $\bU_j$ are nodal values, $N_j$ are polynomial shape functions of order $p$ in the reference element and $\nen$ is the number of nodes per element. The so-called isoparametric mapping, given by
\begin{equation}\label{eqIsoparametricMapping}
\begin{aligned}
\isomap :
\refElem \subset \Rd \ &\longrightarrow \Omega_e \subset \Rd \\
\bxi &\longmapsto
\isomap (\bxi):= \sum_{j=1}^{\nen} \bx^e_j N_j(\bxi),
\end{aligned}
\end{equation}
is employed to establish the relation between the reference element, $\refElem$, and a generic physical element, $\Omega_e$, with nodes $\{\bx_j\}_{j=1,\ldots,\nen}$.
Employing the isoparametric mapping of Eq.~\eqref{eqIsoparametricMapping}, the element and boundary integrals are mapped to the reference space. By using the approximation of the displacement field given by Eq.~\eqref{eqApprox} and selecting the space of weighting functions to be equal to the space spanned by the interpolation functions, the following system of ordinary differential equations is obtained
\begin{equation} \label{eqSystemFEM}
\bM \ddot{\bU} + \bK \bU = \bF.
\end{equation}
As usual in a finite element context, the global mass matrix $\bM$, the global stiffness matrix $\bK$ and the global forcing vector $\bF$ are obtained by assembling the elemental contributions given by
\begin{equation}\label{eqElemContributions}
\bM^e = \int_{\refElem} \rho \bN^T \bN |\isoJ| d\hat{\Omega},
\qquad
\bK^e = \int_{\refElem} (\bB^e)^T \bD^e \bB^e |\isoJ| d\hat{\Omega},
\qquad
\bF^e = \int_{\refElem} \bN^T \bb  |\isoJ| d\hat{\Omega} + \int_{\hat{\Gamma}} \bN^T \bt \| \isoJF \|d\hat{\Gamma}.
\end{equation}
In the above expressions $\bB^e := (\isoJ)^{-1}\grad_S \bN$ is the strain-displacement matrix, $\isoJ$ is the Jacobian of the isoparametric mapping, $\isoJF$ is the Jacobian of the restriction of the isoparametric mapping to an element face and the matrix $\bN$, in three dimensions, is given by
\begin{equation}
\bN := \begin{bmatrix}
N_1 & 0   & 0   & N_2 & 0   & 0   & \cdots & N_{\nen} & 0        & 0 		  \\
0   & N_1 & 0   & 0   & N_2 & 0   & \cdots & 0        & N_{\nen} & 0 		  \\
0   & 0   & N_1 & 0   & 0   & N_2 & \cdots & 0        & 0        & N_{\nen} 
\end{bmatrix}.
\end{equation}
%
\section{The inertia relief method}\label{secIR}

In this section, a short review of the IR method is presented. As already mentioned, the IR method~\cite{wijker2004mechanical} is available in many commercial FE packages and it is widely used in industry to solve unconstrained structural problems without resorting to the more expensive full dynamic analysis. 

When constant unbalanced external loads are applied to an unconstrained structure, the whole system undergoes a steady-state rigid body acceleration in each free direction and, due to the mass of the system, inertial forces are generated that deform elastically the body. The trajectory as rigid body of the system (as if it was infinitely stiff) is described by a displacement field $\bU_r(t)$ such that $\bK \bU_{r}(t)=\boldsymbol{0}$, and therefore also $\bK\dot\bU_{r}(t)=\bK\ddot\bU_{r}(t)=\boldsymbol{0}$. 
The global displacement $\bU$ has to be complemented with the elastic deformation, namely
\begin{equation} \label{eqSplitU}
\bU = \bU_r + \bU_e.
\end{equation}
The elastic deformation field $\bU_{e}$ is important to analyze the internal stresses created by the motion and also to assess other quantities of interest like the torsional stiffness, which is one of the aims of this paper.

The key idea of the IR method is to compute $\bU_{e}$ solving a {\textit{static}} problem 
\begin{equation} \label{eqEquilProb1}
\bK \bU_{e} = \bF_{\text{eq}}
\end{equation}
where the forces $\bF_{\text{eq}}$ are equilibrated, that is the resultant forces and moments are zero. 
Despite matrix $\bK$ is singular, the fact that $\bF_{\text{eq}}$ is equilibrated guaranties that system \eqref{eqEquilProb1} is solvable with a family of infinite solutions, all equivalent up to a rigid body motion. Isostatic constrains (as many as rigid body modes, 3 in 2D and 6 in 3D) have to be set to compute one of these solutions (they all produce the same strains and stresses).

The idea of the inertia relief is to compute the equilibrated forces as
\begin{equation} \label{eqEquilForce1}
\bF_{\text{eq}} = \bF - \bM \ddot\bU_{r}
\end{equation}
noting that $\ddot\bU_{r}$ is the rigid body mode (recall that $\bK\ddot\bU_{r}(t)=\boldsymbol{0}$) such that $\bF_{\text{eq}}$ is equilibrated. Thus, equation \eqref{eqEquilProb1} is derived from \eqref{eqSystemFEM} assuming that $\ddot\bU_{e}=\boldsymbol{0}$ (which stands under a constant load, and therefore constant acceleration).

%
%

\begin{figure}[h]
	\centering
	\includegraphics[width=0.6\textwidth]{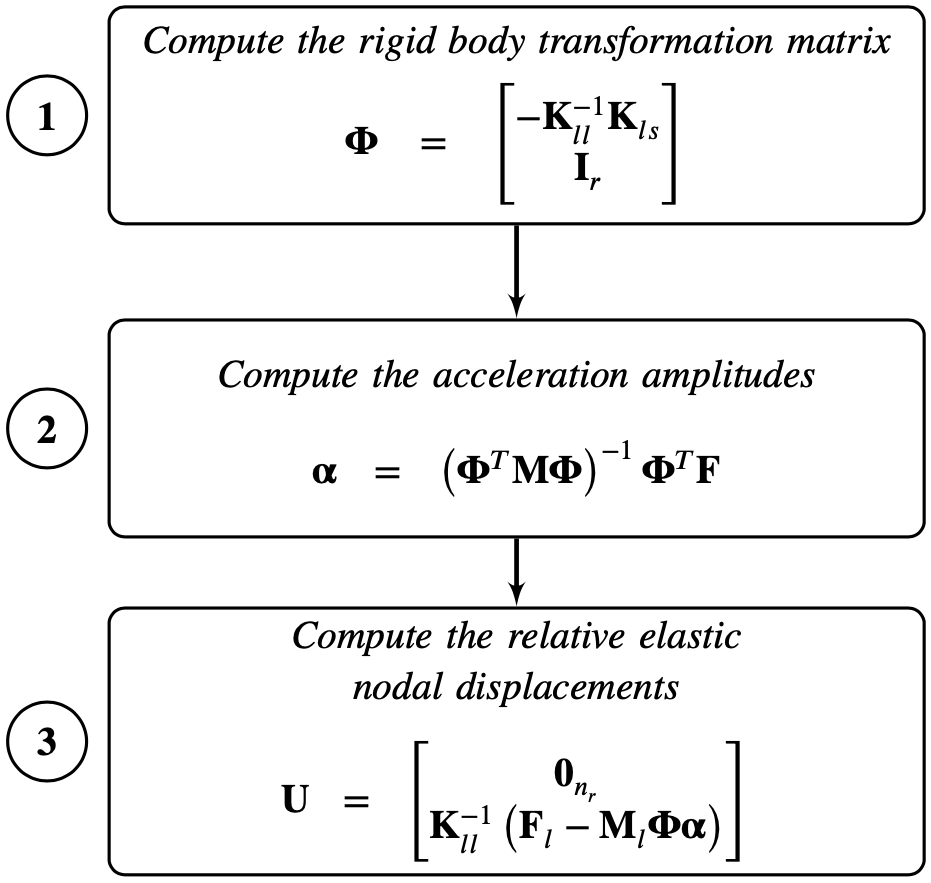}
	\caption{Schematic representation of the steps performed in the IR method.
		\label{InRelSteps}}
\end{figure}

The rigid body acceleration vector $\ddot{\bU}_r$ can be expressed as a linear combination of the (6 in 3D) rigid body modes, namely
\begin{equation}\label{eqAccelR}
\ddot{\bU}_r = \bPhi  \balpha,
\end{equation} 
where each column \textit{rigid body transformation matrix} $\bPhi$ corresponds to one of the $n_r$ rigid body modes ($n_r = 6$ in three dimensions) and the coefficient vector $\balpha$ is seen as containing the acceleration of each of the rigid body modes.
%
Introducing the expression of Eq.~\eqref{eqAccelR} in Eqs.~\eqref{eqEquilForce1} and \eqref{eqEquilProb1}, and pre-multiplying by $\bPhi^T$, the following equation is obtained
\begin{equation} \label{eqStaticIR2}
\bPhi^T \bF - \bPhi^T \bM \bPhi  \balpha = \mathbf{0},
\end{equation}
which guaranties that the right-hand side term in Eq.~\ref{eqEquilProb1} is an equilibrated system of forces (sum of forces and sum of moments equal to zero). It is worth noting that $\bPhi^T \bK = \mathbf{0}$ because the eigenmodes are mutually orthogonal and the eigenvalues (frequencies) associated to the rigid body modes are zero.
The vector of unknown accelerations $\balpha$ providing the equilibrated forces  $\bF_{\text{eq}}$ is therefore computed by solving the system
\begin{equation} \label{eqAccelRef}
\balpha = \left( \bPhi^T \bM \bPhi \right)^{-1} \bPhi^T \bF,
\end{equation}
where $ \bPhi^T \bM \bPhi $ and $\bPhi^T \bF$ are a reduced $6\times 6$ mass matrix and a reduced $6\times 1$ load vector, respectively.
To completely define the rigid body acceleration vector $\ddot{\bU}_r$ in Eq.~\eqref{eqAccelR} it is only necessary to compute the rigid body modes of the structure, that is the 6 columns of matrix $\bPhi$. They correspond to the kernel of the global stiffness matrix, so they are computed as the solution of $\bK \bPhi = \mathbf{0}$. To this end, the set of indices corresponding to the $\ndof = d\times\nN$ degrees of freedom, with $\nN$ being the number of mesh nodes, is partitioned into the reference set $s$ and the remaining set $l$. To simplify the notation, and without loss of generality, the set $s$ is assumed to correspond to the last $\nr$ degrees of freedom. The system of equations to obtain the rigid body modes is then written as
\begin{equation} \label{eqRigidModes}
\begin{bmatrix}
\bKll &  \bKls \\
\bKsl &  \bKss
\end{bmatrix}    
\begin{bmatrix}
\bPhi_{l} \\
\bPhi_{s}
\end{bmatrix}       
= 
\begin{bmatrix}
\mathbf{0}_l  \\
\mathbf{0}_s 
\end{bmatrix}.
\end{equation}
As $\bKll$ is symmetric and positive definite, the degrees of freedom of the rigid body modes corresponding to the $l$ set can be expressed in terms of the degrees of freedom of the rigid body modes corresponding to the reference set, namely
\begin{equation} \label{eqRigidModesSol}
\bPhi_{l} = -\bKll^{-1} \bKls \bPhi_{s}.
\end{equation}
A natural assumption consists of choosing $\bPhi_{s} = \bI_{\nr}$, where $\bI_{\nr}$ denotes the identity matrix of dimension $\nr \times \nr$, so that each column of the matrix $\bPhi_{s}$ represents a unit translation or rotation in the direction of the corresponding reference degrees of freedom.
With all these premises, the relative elastic displacement $\bU_e$ in Eq.~\eqref{eqEquilProb1} is computed. It is worth noting that, in the IR framework, displacements are measured relative to the moving reference set of degrees of freedom \textit{s}, which is subjected to a constant acceleration and undergoes infinite displacements. As a consequence, the rigid body displacement $\bU_r$ is not of interest and can be eliminated from the solving equation. 
Finally, the system to be solved to compute the relative elastic displacement, which in the remainder is simply referred to as $\bU$, reads
\begin{equation} \label{eqSystemIR}
\begin{bmatrix}
\bKll &  \bKls \\
\bKsl &  \bKss
\end{bmatrix}    
\begin{Bmatrix}
\bU_{l} \\
\bU_{s}
\end{Bmatrix}       
= 
\begin{Bmatrix}
\bF_l \\
\bF_s
\end{Bmatrix}  
-
\begin{bmatrix}
\bMll &  \bMls \\
\bMsl &  \bMss
\end{bmatrix} 
\bPhi
\balpha.
\end{equation}
Imposing a zero displacement in the degrees of freedom of the $s$ set, $\bU_{s}=\mathbf{0}$, ensures that the following system is solvable and provides the required relative elastic displacement at the degrees of freedom of the $l$ set,
\begin{equation} \label{eqSystemIR2}
\bKll  \bU_{l} = \bF_l - \bM_l \bPhi \balpha,
\end{equation}
where $\bM_l = \left[ \bMll \;  \bMls\right]$.
%
The IR method can be summarized in three steps as depicted in Fig.~\ref{InRelSteps}. 
%
\section{The parametric inertia relief method}\label{secParamIR}
%
\subsection{Problem definition}\label{secParamIRproblem}

The IR problem is now extended to the case of an unconstrained structure characterised by parametric properties. Let us introduce a set of $\np$ material or geometric parameters denoted by $\bmu = [\mu_1, \mu_2, \dots, \mu_{\np}]^T \in \bmM \subset \mathbb{R}^{\np}$. The set of parametric domains $\bmM$ is defined as the Cartesian product of a predefined interval for each one of the parameters, namely $\bmM := \mM_1\times\mM_2\times\dotsb\times\mM_{\np}$, with $\mu_j\in\mM_j$ for $j=1,\dotsc , \np$ . 
The semi-discrete system of Eq.~\eqref{eqSystemFEM} for a parametric problem can be written as
\begin{equation} \label{eqSystemFEMparam}
\bM(\bmu) \: \ddot{\bU}(\bmu) + \bK(\bmu) \: \bU(\bmu) = \bF(\bmu).
\end{equation}

In order to solve Eq.~\eqref{eqSystemFEMparam} with the IR method, three parametric steps have to be performed, following the rationale of the IR method described in Sec.~\ref{secIR} for the non-parametric case. The first step consists of computing the rigid body modes as
\begin{equation} \label{eqIRparam1}
\bPhi(\bmu) = 	
\begin{bmatrix}
- \bK^{-1}_{ll}(\bmu) \: \bK_{ls}(\bmu)  \\
\bI_r
\end{bmatrix}.
\end{equation}
Second, the vector of accelerations is given by 
\begin{equation} \label{eqIRparam2}
\balpha(\bmu) = \left[ { \bPhi}^{T}(\bmu) \: \bM(\bmu) \: \bPhi(\bmu) \right]^{-1} \: \bPhi^T(\bmu) \: \bF(\bmu).
\end{equation}
Finally, the relative elastic displacement is computed as
\begin{equation} \label{eqIRparam3}
\bU (\bmu)  = 
\begin{bmatrix}
\mathbf{0}_{\nr} \\
\bKll^{-1}(\bmu)  \left( \bF_l(\bmu) - \bM_l(\bmu)  \: \bPhi(\bmu) \: \balpha(\bmu)  \right)
\end{bmatrix}.
\end{equation}
In Eqs.~\eqref{eqIRparam1} to~\eqref{eqIRparam3}, $\bmu$ is treated as a set of additional independent variables (or parametric coordinates), instead of problem parameters. As a consequence, the generalized solution of the three equations depends explicitly on the parameters and takes values in the multidimensional domain $\mathcal{D} = \Omega \times \bmM$. Standard numerical methods (e.g. finite elements, finite volumes, finite differences) would require the solution of each step of the IR method in the high dimensional space $\mathcal{D}$, which is not feasible in practical problems. In this work, the PGD is proposed as a ROM able to circumvent the so-called \textit{curse of dimensionality} and to provide the generalized solution of the parametric IR problem. 

\subsection{Cascade application of the encapsulated PGD approach}\label{secParamIRpgd}

The goal of this section is to solve the parametric IR problem by means of the PGD technique. Following the standard PGD rationale, let us assume that the solution $\fU(\bmu)$ of  Eq.~\eqref{eqSystemFEMparam} can be approximated by a linear combination of an \textit{a priori} unknown number $N_{\bU}$ of terms (or modes), namely
\begin{equation}\label{eqUPGD}
\fU(\bmu) \approx \fU^{\texttt{PGD}}(\bmu) =\sum_{i=1}^{N_{\bU}} \beta_{\bU}^i \: \fU^i \: \prod_{j=1}^{n_p} u_{ j}^i(\mu_j).
\end{equation} 
Each PGD mode $i$ is given by the product of a spatial term, $\fU^i$, defined on the discretized space $\Omega$ and a set of parametric functions $u_{ j}^i(\mu_j)$ depending, in a separated form, on each parameter $\mu_j$, for $j = 1,2,\dots, \np$. The spatial term, $\fU^i$, is a vector of the size of the finite element displacement vector, whereas each parametric dimension $\mu_j$ is discretized with $n_j$ points with coordinates $\mu_j^{p_j}$, where $p_j = 1, 2, \dots, n_j$. The spatial and parametric modes are usually normalized and the amplitude of each mode, $\beta_{\bU}^i$, indicates the relevance of the $i$-th mode to the final separated solution.

In order to compute the terms of the summation in Eq.~\eqref{eqUPGD}, the PGD solver typically employs a greedy approach. Assuming that the previous $n-1$ modes are known, the greedy algorithm computes sequentially the $n$-th term 
\begin{equation}\label{eqGreedy}
\fU^{\texttt{PGD},n}(\bmu) = \fU^{\texttt{PGD},n-1}(\bmu)  +  \fU^n\: \prod_{j=1}^{n_p} u^n_{j}(\mu_j)
\end{equation} 
given by the spatial mode  $\fU^n$ and the parametric terms $u_{j}^n(\mu_j)$ for $j = 1,2,\dots, \np$.
The enrichment process automatically stops when a user-defined level of accuracy is reached, that is when the amplitude $\beta_{\bU}^n$ of the last term is smaller than a user defined tolerance. 

Since the unknown spatial and parametric terms $\fU^n$ and $u^n_{j}(\mu_j)$ are multiplying, the problem of computing the $n$-th term in Eq.~\eqref{eqGreedy} is nonlinear. More precisely, it is a nonlinear least-squares problem defined to find the best rank-one approximation (meant as the product of sectional functions) of the unknown term $\fU^n\: \prod_{j=1}^{n_p} u^n_{j}(\mu_j)$. As commonly done in a PGD context, an alternated direction scheme is applied, which consists in solving the problem separately for each unknown function, assuming that all the others are known, until a stationary solution is reached. It is worth emphasizing that despite a non-linear problem needs to be solved to obtain each PGD mode, the computational cost of the problem increases linearly with the number of introduced parameters, making the solution of high-dimensional problems affordable. 
Recently, D\'{i}ez et al. \cite{diez2019encapsulated} developed the \textit{Encapsulated PGD Toolbox}, which is a collection of PGD-based algorithms able to perform algebraic operations (e.g. product, division, storage, compression, solving linear system of equations, etc.) for multidimensional data represented in a discretized tensorial separated format. The main advantage of the library (freely available at \url{https://git.lacan.upc.edu/zlotnik/algebraicPGDtools.git}) is that each routine is \textit{encapsulated}, meaning that it can be used as a black box. This is particularly attractive for the end user and it facilitates the interaction with commercial software. 

\begin{figure}[h]
	\centering
	\includegraphics[width=1.\textwidth]{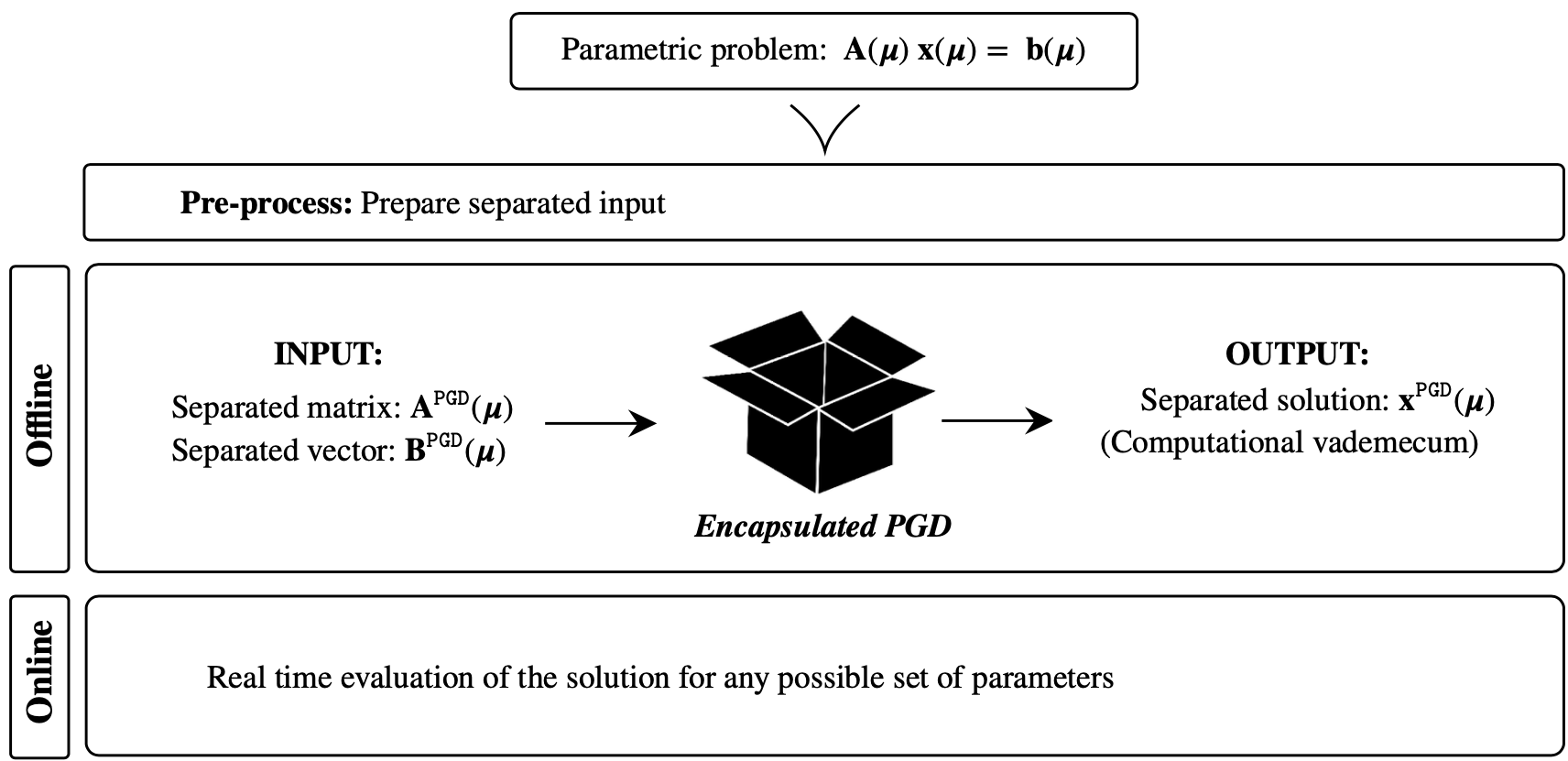}
	\caption{Structure of the \textit{Encapsulated-PGD} linear solver. \label{fig:EncapsPGD}}
\end{figure}

\begin{figure}[h]
	\centering
	\includegraphics[width=\textwidth]{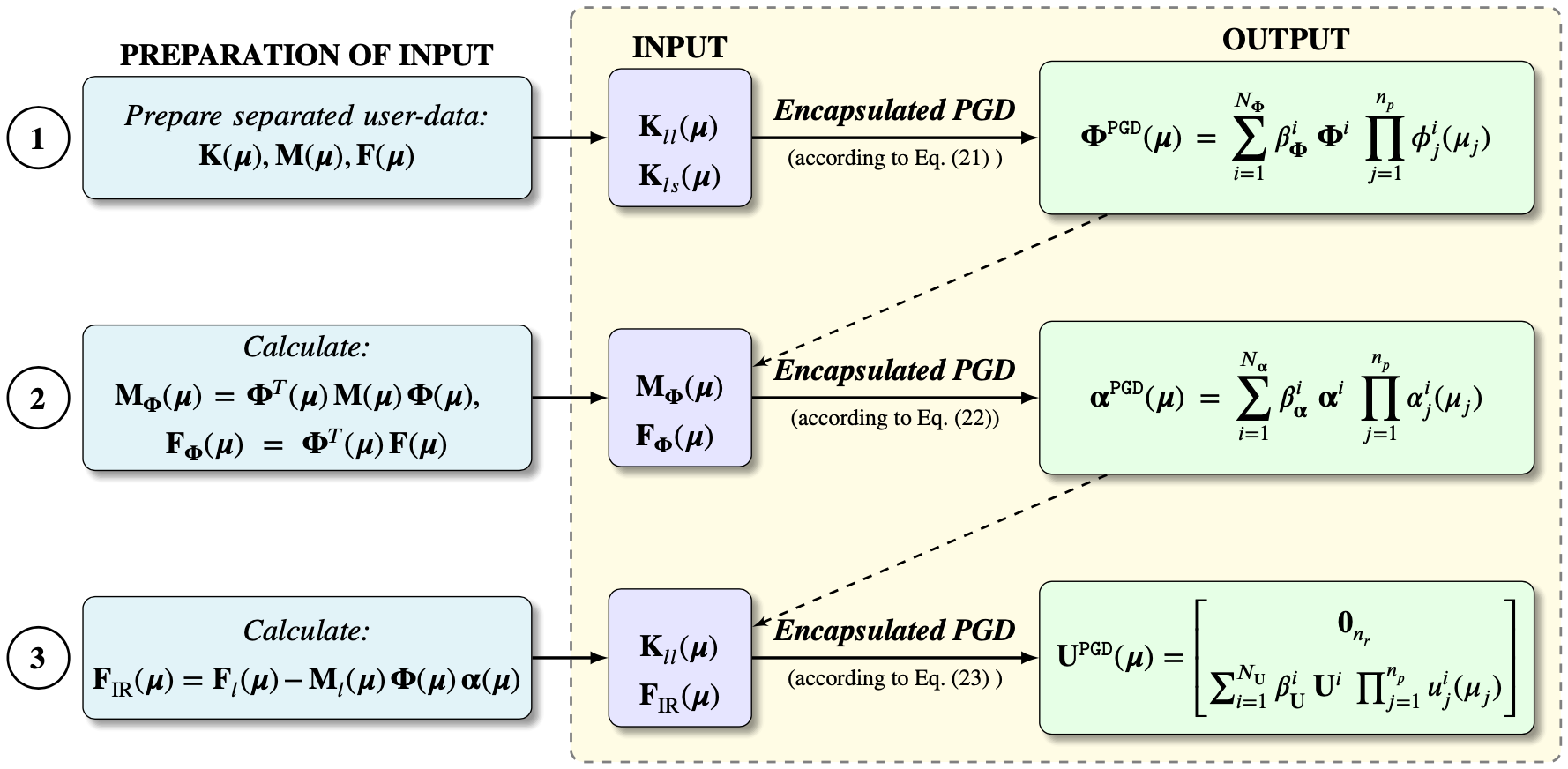}
	\caption{Schematic representation of the \textit{cascade} encapsulated PGD approach for the solution of a parametric IR problem. }\label{fig:ConcatPGD}
\end{figure}

To illustrate the idea behind the encapsulated PGD Toolbox, Fig.~\ref{fig:EncapsPGD} describes the structure of the encapsulated-PGD routine that solves parametric linear systems of equations. In a straightforward way, the same structure can be extended to other arithmetic operators. Shortly, given an algebraic linear system of equations $\fA(\bmu)  \: \fx(\bmu) = \fb(\bmu)$ depending on the set of parameters $\bmu$, the toolbox is able to return an explicit description of $\fx (\bmu)$, also called \textit{computational vademecum}, containing the solution for every possible combination of the parameters. The only requirement to employ the encapsulated PGD approach is to pre-process the input quantities, i.e. the parametric matrix $\fA (\bmu) $ and vector $\fb (\bmu)$, such that they are expressed in a PGD separated form. Given the input data, the user only needs to employ the encapsulated PGD in the offline stage, to obtain the PGD approximation by means of the above mentioned greedy algorithm and alternate direction scheme. The output consists of  the sought \textit{computational vademecum} and, during an online stage, the solution can be evaluated in real time for any set of parameters at a negligible computational cost.
The PGD-IR approach proposed in this work makes use of the encapsulated PGD toolbox. The three parametric IR Eqs.~\eqref{eqIRparam1} to~\eqref{eqIRparam3} are solved sequentially, in the sense that the solution of each equation is directly needed in the next one. Consequently, as depicted in Fig.~\ref{fig:ConcatPGD}, a \textit{cascade} PGD scheme can be employed, in which the output of each step, obtained in a separated format by simply calling the encapsulated PGD linear solver, can be directly used as an input of the next one, until the final solution of the global problem is computed. 
It is worth noting that, in order to use the toolbox, the user has to provide a separable representation of the input data. In particular, the stiffness $\fK(\bmu)$ and mass $\fM(\bmu)$ matrices must be written as
\begin{align}\label{eqKMseparated}
\fK(\bmu)&\approx \sum_{i=1}^{N_{\fK}} \fK^i  \prod_{j=1}^{n_p} k_{j}^i(\mu_j), \\
\fM(\bmu)&\approx \sum_{i=1}^{N_{\fM}} \fM^i  \prod_{j=1}^{n_p} m_{j}^i(\mu_j),  
\end{align}
where the spatial terms are $\fK^i \in \mathbb{R}^{\ndof \times \ndof}$ and $\fM^i \in \mathbb{R}^{\ndof \times \ndof}$, whereas the parametric terms are $k_{j}^i(\mu_j) \in \mathbb{R}^{n_j}$ and $m_{j}^i(\mu_j) \in \mathbb{R}^{n_j}$, for $j = 1,2, \dots, n_p$. Similarly, the input nodal force vector $\fF(\bmu)$ must be written as
\begin{equation}\label{eqFseparated}
\fF(\bmu)\approx \sum_{i=1}^{N_{\fF}} \fF^i  \prod_{j=1}^{n_p} f_{j}^i(\mu_j),
\end{equation}
with $\fF^i \in \mathbb{R}^{\ndof }$ and $f_{j}^i(\mu_j) \in \mathbb{R}^{n_j}$. In the above expressions $N_{\fK}$, $N_{\fM}$ and $N_{\fF}$ are the number of modes required to produce a separable approximation of $\fK(\bmu)$, $\fM(\bmu)$ and $\fF(\bmu)$ respectively.

It is important to underline that it is not always trivial to find a separated representation of the input data, as given by equations~\eqref{eqKMseparated} and ~\eqref{eqFseparated}, especially when geometric parameters are considered in the problem. This issue will be addressed in the next section.

Finally, in order to solve the parametric equations of the PGD-IR approach, steps 2 and 3 depicted in Fig.~\ref{fig:ConcatPGD} require some extra operations between the parametric objects, such as products, additions or compression. These operations can be easily performed by the \textit{Encapsulated PGD toolbox}. 
\begin{remark} \label{rk:materialPGDIR}
	The first two steps of the PGD-IR shown in Fig.~\ref{fig:ConcatPGD} are parametric problems only if geometric parameters are considered, because, by definition, the rigid body modes of a structure do not depend on the material properties. 
\end{remark}

\subsection{Geometric parameters: a nonintrusive algebraic approach to separate input quantities}\label{secAlgSep}

The extension of the proposed nonintrusive PGD framework to geometrically parametrized problems represents a challenging task. This is due to the fact that, if geometric parameters are introduced in the problem, it is not trivial to find separable representation of the input quantities. 

If a closed form separated expression of the stiffness and mass matrices is sought, the weak form of the problem must be modified to account for the parametric geometry. A common approach consists of formulating the problem in a reference domain, leading to several limitations that are briefly discussed in Appendix~\ref{app1}. The most important limitation in the context of the current work is that the implementation based on a reference domain requires access to the code, precluding its application in an industrial framework, where commercial codes are typically employed.

In this section a nonintrusive algebraic approach is proposed, which is able to deal with general geometric parametrizations. The main idea is to perform a sampling of the parametric matrices and to express them in a separated format. The approach requires the computation of the parametric matrices for different values of the geometric parameters, whilst maintaining the connectivity matrix of the FE mesh. To this end a mesh morphing approach is adopted in this work. Every time a sampling of the parametric matrices is required, an initial mesh is deformed according to the geometric parameters and the global stiffness and mass matrices are computed. It is worth noting that this approach can be easily integrated in commercial packages that are equipped with a mesh morphing capability. Alternatively, the user can define the preferred mesh morphing approach and produce a set of meshes to be imported in the preferred FE software. It is also worth mentioning that the sampling does not require the solution of the FE system of equations as only the global stiffness and mass matrices are of interest for the proposed PGD-IR approach. Once the set of global stiffness and mass matrices is available, they are expressed in a separated format using the encapsulated PGD toolbox.

To illustrate the proposed nonintrusive approach, let us consider the stiffness matrix $\fK\in \mathbb{R}^{\ndof \times \ndof}$, depending on $n_p$ parameters $\bmu = [ \mu_1, \mu_2, \dots, \mu_{\np}]^T \in \bmM \subset \mathbb{R}^{n_p}$. The parametric dimension $\mu_j \in \mM_j$, for $j = 1,2,\dots, \np$, is discretized using $n_j$ points with coordinates $\mu_j^{p_j}$, where $p_j = 1,2,\dots,n_j$. The full-order sampling of the parametric matrix consists of evaluating $\fK (\bmu)$ in the set of $n_{\text{tot}}$ points used to discretize the parametric domain $\bmM = \mM_1\times\mM_2\times\dotsb\times\mM_{\np}$, where $n_{\text{tot}} = \prod_{j=1}^{n_p} n_j$. Each point is characterized by its sectional indices $(p_1,p_2,\dots,p_{\np})$, which are duly sorted by using a linear array index $i$ such that
\begin{equation}\label{eqIndex}
i = p_1 + (p_2 -1) n_2 + (p_3 -1) n_2 \times n_3  + \dots = p_1 + \sum_{j=2}^{n_p} (p_j -1) \prod_{l = 2}^{j} n_l.
\end{equation}

Note that the association between the multi index $(p_1,p_2,\dots,p_{\np})$ and the index $i$ is also obtained by updating $i = i+1$ inside $n_p$ nested loops, with no need to use explicitly expression~\eqref{eqIndex}. Employing the association between the multi-index $(p_1,p_2,\dots,p_{\np})$ and the linear index $i$, the parametric stiffness matrix $\bK(\bmu)$ can be written as
\begin{equation}\label{eqKassociation}
\bK(\bmu) = \sum_{p_1 =1}^{n_1} \sum_{p_2 =1}^{n_2} \dots \sum_{p_{n_p} =1}^{n_{n_p}} \fK (\mu_1^{p_1}, \mu_2^{p_2}, \dots , \mu_{n_p}^{p_{n_p}}) \: F_{p_1 , p_2, \dots, p_{n_p}} (\mu_1, \mu_2, \dots, \mu_{\np}),
\end{equation}
where $F_{p_1 , p_2, \dots, p_{n_p}}$ is such that $F_{p_1 , p_2, \dots, p_{n_p}} (\mu_1^{p_1}, \mu_2^{p_2}, \dots, \mu_{\np}^{p_{n_p}}) = 1$ and it is equal to zero for any other values of the discrete indices $p_j$. 
Using the linear indexing $i$ introduced in Eq.~\eqref{eqIndex}, Eq.~\eqref{eqKassociation} becomes
\begin{equation}\label{eqKseparatedAlgebraic}
\fK(\bmu) = \sum_{i=1}^{n_{\text {tot}}} \fK^i \prod_{j=1}^{n_p} k_j^i(\mu_j),
\end{equation}
where $\fK^i = \fK(\mu_1^{p_1}, \mu_2^{p_2}, \dots , \mu_{n_p}^{p_{n_p}})$, and $k_j^i (\mu_j^{p_l}) = \delta_{p_l,p_j}$ for any $p_l = 1,2, \dots, n_j$, while $p_j$ is given by $i$ as defined in Eq.~\eqref{eqIndex}. Finally, Eq.~\ref{eqKseparatedAlgebraic} represents the desired separated representation of the stiffness matrix. 

Depending on the number of parameters, $n_p$, and the number of nodes chosen to discretize the parametric domains, $n_j$, the separated expression of the parametric stiffness matrix might involve a large number of terms, $n_{\text{tot}}$. It is possible to reduce the computational cost of the following calculations by employing the PGD-compression, available in the encapsulated PGD-toolbox~\cite{diez2019encapsulated}. The idea is to perform an $\mathcal{L}_2$ projection of the expression of Eq.~\eqref{eqKseparatedAlgebraic} to reduce the number of terms in the summation while maintaining an accurate representation of $\fK(\bmu)$.
In a similar fashion, a separated representation of the mass matrix can be also obtained. As it will be shown by means of numerical examples, the main advantage of the proposed algebraic technique is its flexibility which in general allows to add an arbitrary number of geometric parameters as variables of the problem. In addition, the nonintrusive character of the proposed ROM, makes the approach proposed in this work suitable for industrial applications.

\section{Numerical examples}\label{secNumEx}

In this section two numerical examples are presented in order to show the properties of the proposed method. The first example is used to illustrate the numerical properties of the proposed PGD-IR method when both material and geometric parameters are considered. In the second example, the method is applied to a more realistic industrial case involving three parameters. Furthermore, a multiobjective optimization study is performed, which proves the potential of the PGD-IR method in the context of design optimization problems.

\subsection{Parametric inertia relief with material and geometric parameters}\label{secNumEx1}

A pure torsion test case is considered for an unconstrained linear elastic 3D structure characterized by one material and one geometric parameter, that are treated as additional coordinates of the problem. For a better readability, the two variables are denoted here with different symbols, that is $\mu \in \mM_{\mu}$ and $\theta \in \mM_{\theta}$ for the material and geometric parameters respectively. 

\begin{figure}[tb!]
	\begin{subfigure}[b]{.55\linewidth}
		\centering
\textsl{}		\includegraphics[width=0.8\linewidth]{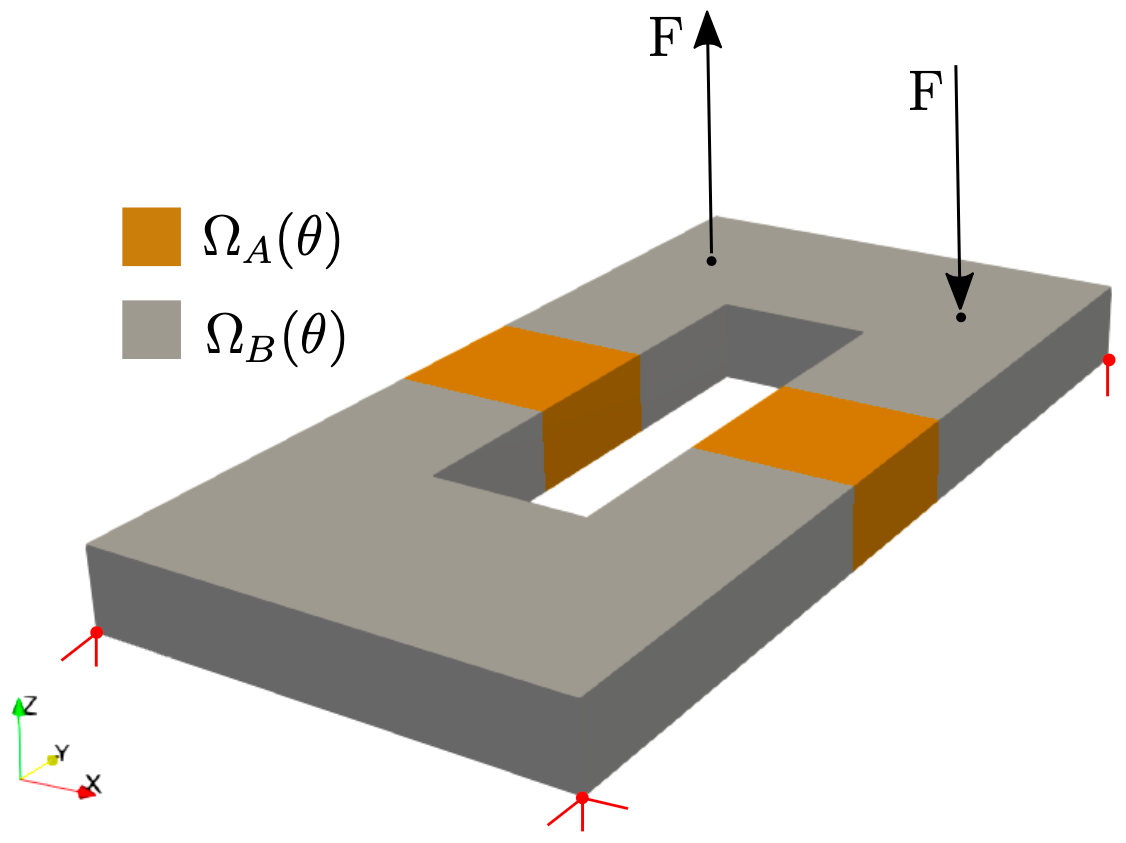}
		\label{fig:donut_set}
	\end{subfigure}%
	\begin{subfigure}[b]{.5\linewidth}
		\centering
		\includegraphics[width=0.8\linewidth]{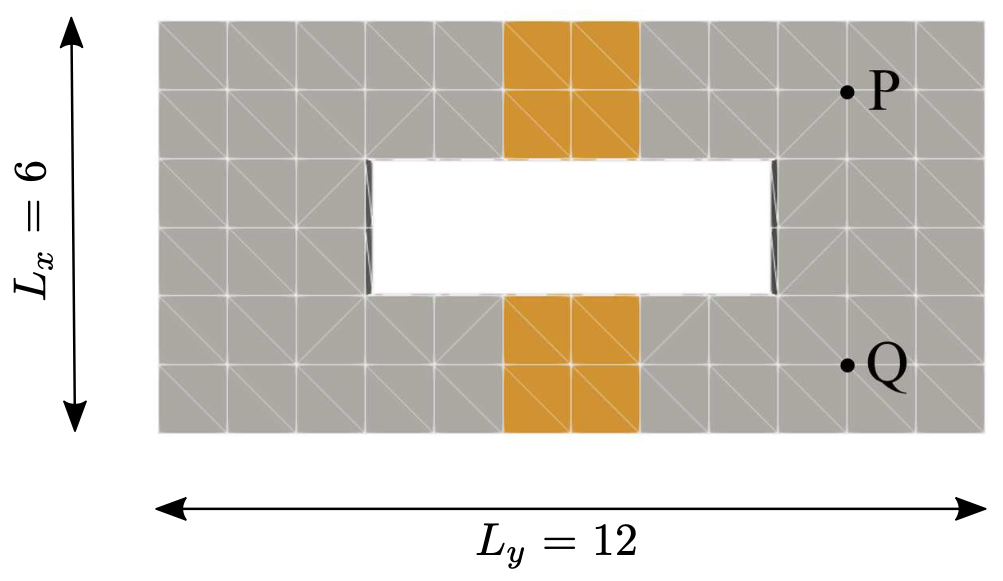}
		\label{fig:donut_mesh}
	\end{subfigure}
	\caption{Computational domain, showing the partition into two non-overlapping subdomains $\Omega_A(\theta)$ and $\Omega_B(\theta)$ (left) and top view of the discretized computational domain showing the dimensions and the points $P$ and $Q$ where the forces are applied (right).}
	\label{fig:Donut}
\end{figure}

As depicted in Fig.~\ref{fig:Donut}, the reference domain $\hat \Omega$ consists of a block with dimensions $[-L_x/2,L_x/2] \times [-L_y/2,L_y/2] \times [-L_z/2,L_z/2]$ with an inclusion given by $[-L_x/6,L_x/6] \times [-L_y/4,L_y/4] \times [-L_z/2,L_z/2]$, where $L_x = 6$, $L_y = 12$ and $L_z = 1$. The torsional load is given by two parallel forces of constant magnitude $F = 10$ acting on the positive and negative $z$ direction respectively and applied at the points $P=(2,4,1/2)$ and $Q=(-2,4,1/2)$. Fig.~\ref{fig:Donut} also shows the spatial discretization employed, consisting on a regular mesh with 236 nodes and 742 linear tetrahedral elements. 

The physical domain $\Omega(\theta)$ depends upon the geometric parameter and it is split into two non-overlapping subdomains $\Omega_A(\theta)$ and $\Omega_B(\theta)$. The parametric Young's modulus $E$ is defined as
\begin{equation}\label{eqMu_Ex1}
E(\bx, \mu) =
\begin{cases}
E_A (\mu) = \mu  \quad &\text{for} \quad \bx \in \Omega_A( \theta), \\
E_B = 200 \quad &\text{for} \quad \bx \in \Omega_B(\theta),
\end{cases}
\end{equation}
where the Young modulus $E_A(\mu)$ is considered varying in the range $\mM_{\mu}= [10, 410]$, and $\mM_{\mu}$ is discretized with a uniform distribution of $n_{\mu} = 41$ points. The Poisson's ratio and the density are assumed constant in the whole domain and taken as $\nu = 0.3$  and $\rho=1$ respectively.

The geometrically parametrized domain $\Omega(\theta)$ is described with the Cartesian coordinates $\bx$, and it is defined as the image of the reference domain $\hat \Omega$, with reference coordinates $\hat \bx$, via a geometric mapping $\bm{\Psi}(\hat \bx, \theta)$, namely
\begin{equation}\label{eqMapping_Ex1}
\begin{cases}
\displaystyle
x &= \psi_1(\hat \bx, \theta) = \hat x + \theta \sin{\left(\displaystyle \frac{\pi \hat y}{L_y}\right)} \left(\hat x - \displaystyle \frac{L_x}{2}\right), \\
y &= \psi_2(\hat \bx, \theta) = \hat y, \\
z &= \psi_3(\hat \bx, \theta) = \hat z.
\end{cases}
\end{equation}
The parameter $\theta$ is taken to be in the interval $\mM_{\theta}=[0, 0.5]$, and $\mM_{\theta}$ is discretized with a uniform distribution of $n_{\theta} = 21$ points. 

\begin{figure}[h]
	\small
	\centering
	\includegraphics[width=1\textwidth]{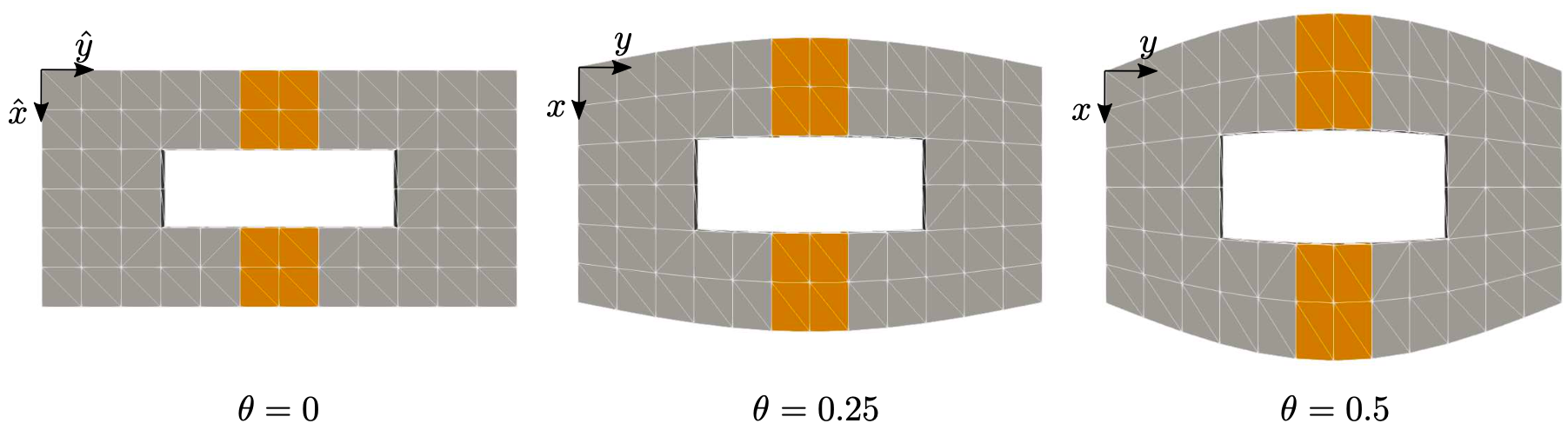}
	\caption{Physical domain for three different values of the geometric parameter $\theta$. }
	\label{fig:donut_deformed}
\end{figure}

Fig.~\ref{fig:donut_deformed} shows the influence of the parameter $\theta$ in the geometry of the computational domain for three different values of $\theta$. The particular value $\theta=0$ leads to a deformed configuration that coincides with the reference configuration, i.e. the mapping of Eq.~\eqref{eqMapping_Ex1} becomes the identity. The configurations in Fig.~\ref{fig:donut_deformed} also show that the mapping changes the nodal coordinates of the mesh while maintaining the connectivities, as required within the current PGD framework.

The objective of this numerical test is to employ the proposed PGD-IR approach to obtain a \textit{computational vademecum} able to describe the variation of the solution with respect to the material and geometric parameters.

Following the proposed PGD-IR framework, the first step consists of choosing a reference set of six degrees of freedom able to counteract the rigid body motions of the structure. Next, in order to employ the encapsulated PGD toolbox, it is necessary to define the input data (i.e. stiffness matrix, mass matrix, force vector) in a separated format. By using the linear dependence of the stiffness matrix on the Young's modulus, an analytical separable representation of the stiffness matrix with respect to $\mu$ can be easily constructed. For the geometric parameter $\theta$, the algebraic PGD toolbox is employed, as discussed in detail in Sec.~\ref{secAlgSep}. For every nodal value of the geometric parameter $\theta^{p} = [\theta^1, \theta^2, \dots, \theta^{n_{\theta}}]^T$, the geometrically deformed mesh is generated according to the mapping of Eq.~\eqref{eqMapping_Ex1}, and two stiffness-like matrices $\fK_A({\theta^{p}})$ and $\fK_B({\theta^{p}})$ are computed. The quantity $\fK_A({\theta^{p}})$ is calculated by imposing the Young's modulus $(E_A, E_B) = (1, 0)$, thus accounting for the contribution of the finite elements belonging to the subdomain $\Omega_A(\theta^{p})$ to the global stiffness matrix. Analogously, $\fK_B{\theta^{p}}$ corresponds to the choice $(E_A, E_B) = (0, 1)$ and accounts for the contribution of the finite elements belonging to the subdomain $\Omega_B(\theta^{p})$. Once these matrices are sampled in the parametric nodes $n_{\theta}$, a separated form of the  parametric global stiffness matrix is readily available, namely
\begin{equation}\label{eqK_ex1}
\fK(\mu,\theta) = E_A(\mu) \sum_{i=1}^{n_{\theta}} \fK_A^i \: k^i(\theta) + E_B \sum_{i=1}^{n_{\theta}} \fK_B^i \: k^i(\theta),
\end{equation}
with $k^i(\theta^{p}) = \delta_{p, i}$, for every $p = 1, 2, \dots, n_{\theta}$. In this example, a PGD-compression was performed, which is always advisable when the number of PGD-terms is large, and an accurate approximation of the stiffness matrix was obtained in the known PGD format
\begin{equation}\label{eqKcomp_ex1}
\fK^{\texttt{PGD}}(\mu,\theta) = \sum_{i=1}^{N_{\fK}} \fK^i \: k^i(\mu) \: k^i(\theta).
\end{equation}
In this case, after performing compression with a tolerance $tol = 10^{-5}$, the number of PGD terms was reduced to $N_{\fK} = 10$.
Following the same procedure, the PGD approximation of the parametric mass matrix is also obtained, namely 
\begin{equation}\label{eqM_ex1}
\fM^{\texttt{PGD}}(\mu,\theta) = \sum_{i=1}^{N_{\fM}} \fM^i \: m^i(\mu) \: m^i(\theta).
\end{equation}

Please note that the mass matrix is actually independent on the Young modulus, that is $m^i(\mu) = 1$. However, the general expression of Eq.~\eqref{eqM_ex1} is used to maintain a consistent notation for all the inputs of the PGD-IR approach. 
Finally, the global forcing vector is also written in the general separated form
\begin{equation}\label{eqF_ex1}
\fF^{\texttt{PGD}} = \fF \: f(\mu) \: f(\theta),
\end{equation}
where, again, it is worth emphasizing that $\fF$ is the standard FE forcing vector and $f(\mu) = f(\theta) = 1$, because the right hand side is not dependent on the material parameter and, for the given set of forces applied to the structure is also independent on the geometric parameter.

The computation of the separated form of the stiffness and mass matrices and the forcing vector completes the pre-process required to apply the proposed PGD-IR approach. Next, the three steps of the PGD-IR approach can be sequentially completed. As detailed in Remark~\ref{rk:materialPGDIR}, the three steps involve a parametric problem because not only material parameters are considered but also geometric parameters, leading to a generalized solution that can be written as
\begin{equation}\label{eqU_PGDex1}
\fU^{\texttt{PGD}}(\mu, \theta) = \sum_{i=1}^{N_{\bU}} \beta_{\bU} \: \fU^i \: u_{\mu}^i(\mu) \: u_{\theta}^i(\theta).
\end{equation}

It is worthy to mention that the proposed PGD-IR approach was implemented in a Matlab routine which acts as a black-box, following the philosophy of the encapsulated PGD toolbox. In fact, the routine only requires to receive the input quantities in a separated form in order to return the output in the same separated form.

\begin{figure}[tb!]
	\centering
	\includegraphics[width=.50\textwidth]{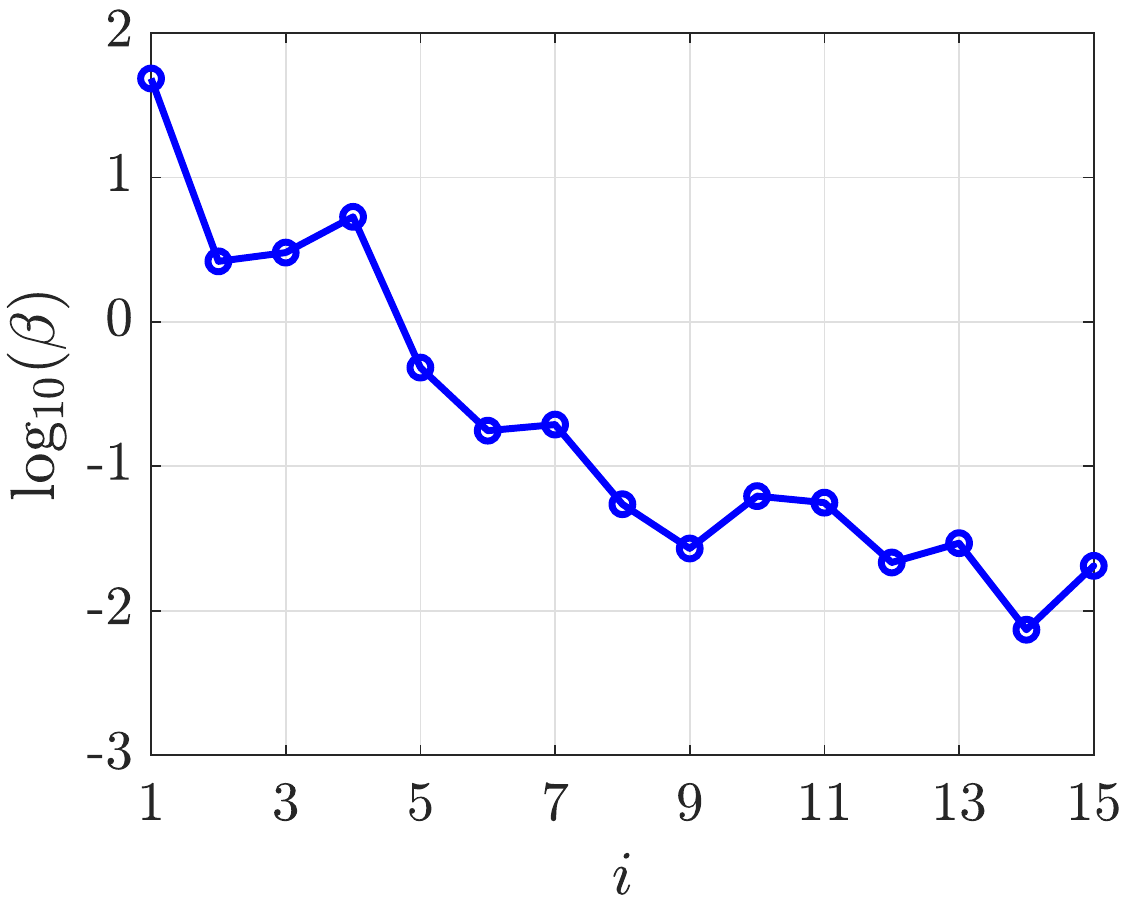}
	\caption{Evolution of the amplitude of the PGD modes $\beta^i$ of the solution $\fU^{\texttt{PGD}}(\mu, \theta)$ with respect to the number of PGD modes, $i$.}
	\label{fig:amplitude}
\end{figure}

Fig.~\ref{fig:amplitude} plots the evolution of the amplitude $\beta_{\bU}$ of each PGD mode. It can be observed that the amplitude rapidly decreases as the number of modes is increased. With 15 computed modes the amplitude of the last mode is almost four orders of magnitude lower than the amplitude of the first mode. In addition, the results show that the first four modes capture the most relevant information of the generalized solution as the fifth and subsequent modes have an amplitude at least two orders of magnitude lower than the amplitude of the first mode. 

\begin{figure}[tb!]
	\begin{subfigure}{.24\linewidth}
		\centering
		\includegraphics[width=1.25\linewidth]{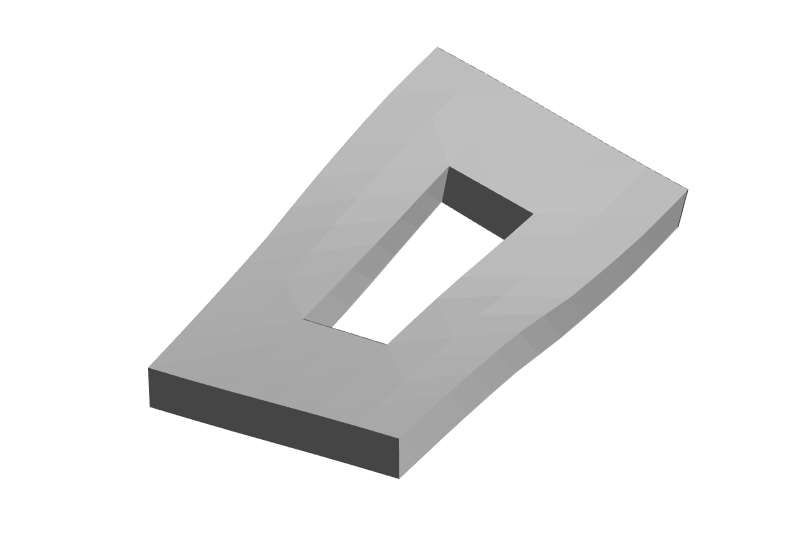}
		\caption{Mode 1}
	\end{subfigure}%
	\begin{subfigure}{.24\linewidth}	
		\centering
		\includegraphics[width=1.25\linewidth]{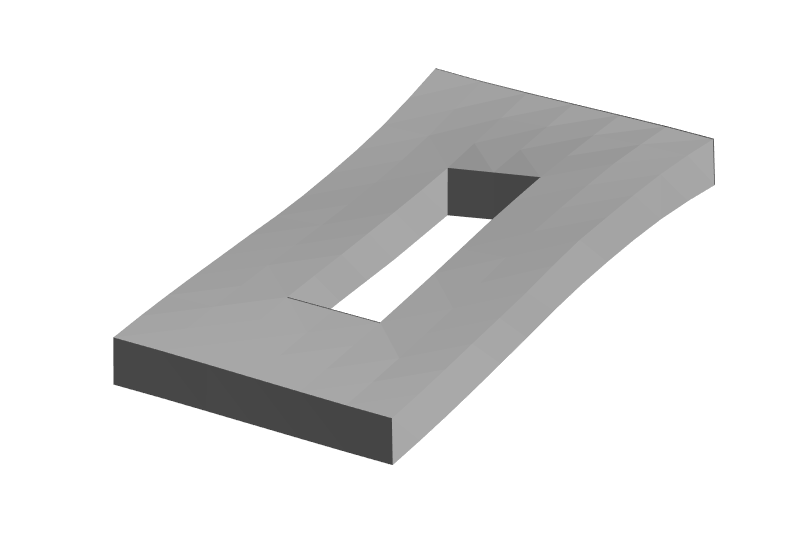}
		\caption{Mode 2}
	\end{subfigure}
	\begin{subfigure}{.24\linewidth}
		\centering
		\includegraphics[width=1.25\linewidth]{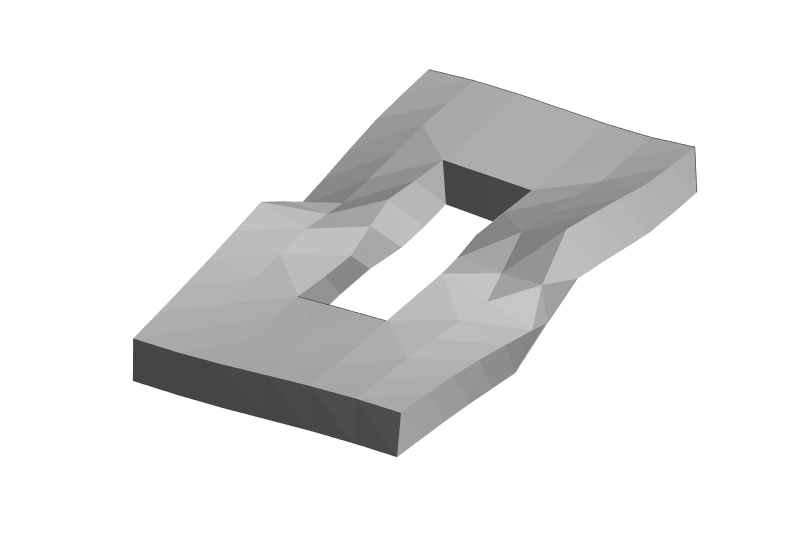}
		\caption{Mode 3}
	\end{subfigure}%
	\begin{subfigure}{.24\linewidth}
		\centering
		\includegraphics[width=1.25\linewidth]{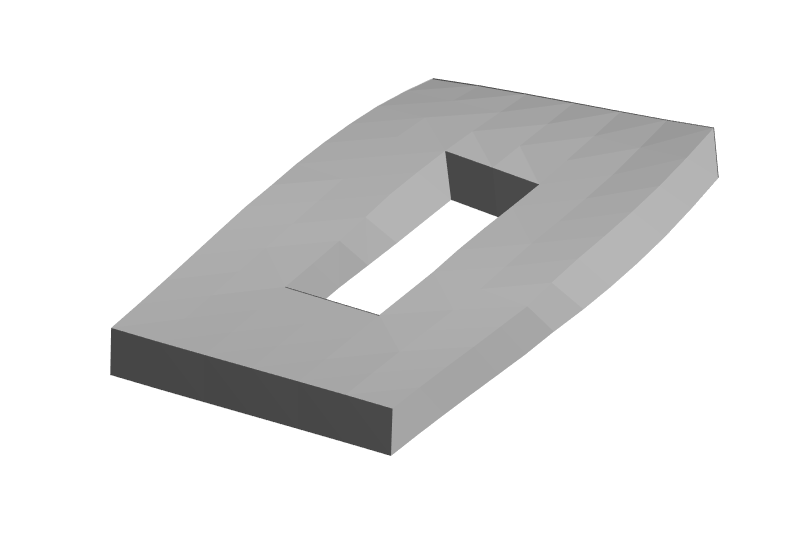}
		\caption{Mode 4}
	\end{subfigure}
	\caption{First four spatial modes of the generalized solution $\fU^{\texttt{PGD}}(\mu, \theta)$.}
	\label{fig:spatialModes}
\end{figure}

\begin{figure}[tb!]
	\begin{subfigure}{.5\linewidth}
		\centering
		\includegraphics[width=1.\linewidth]{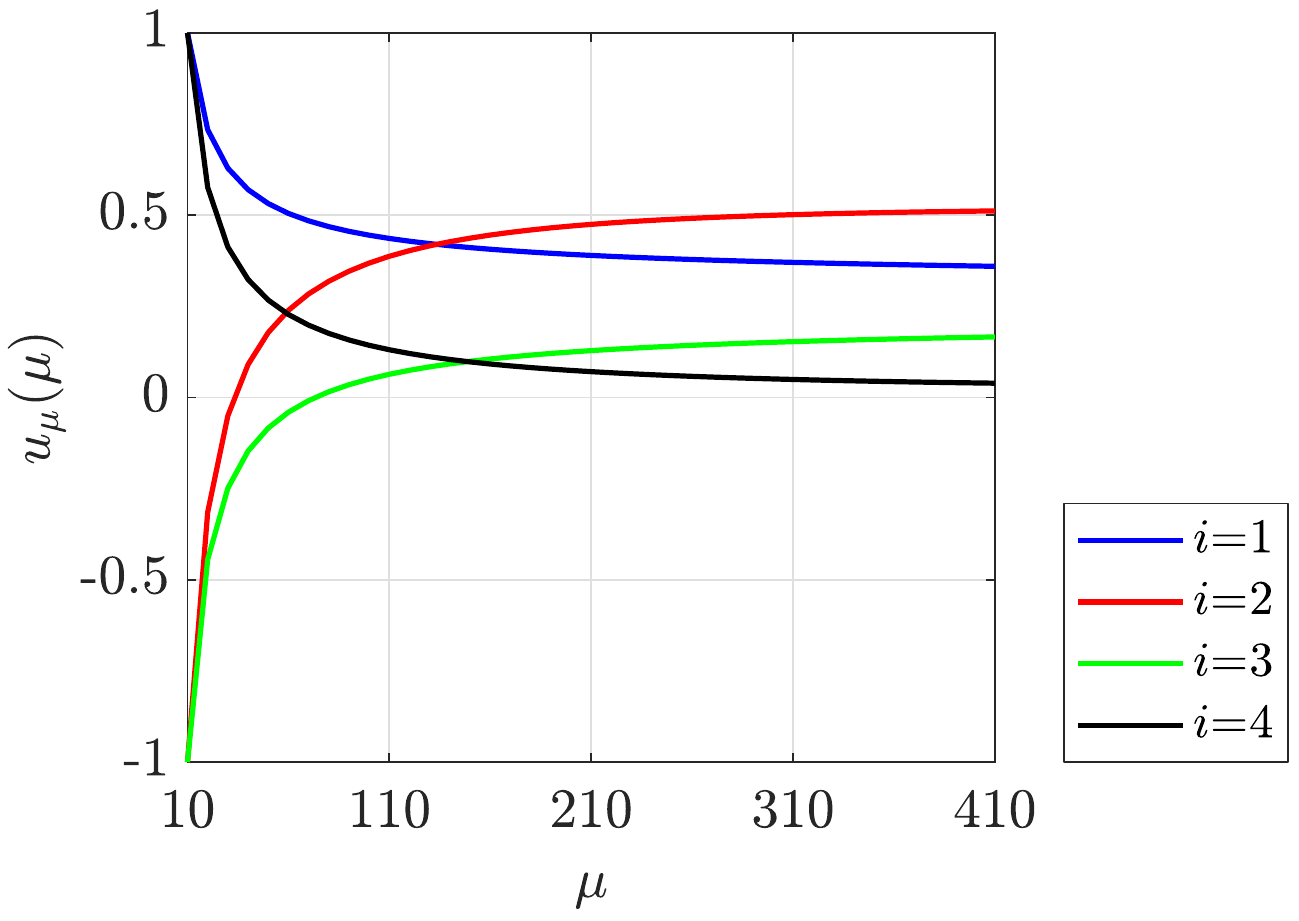}
	\end{subfigure}%
	\begin{subfigure}{.5\linewidth}
		\centering
		\includegraphics[width=1.\linewidth]{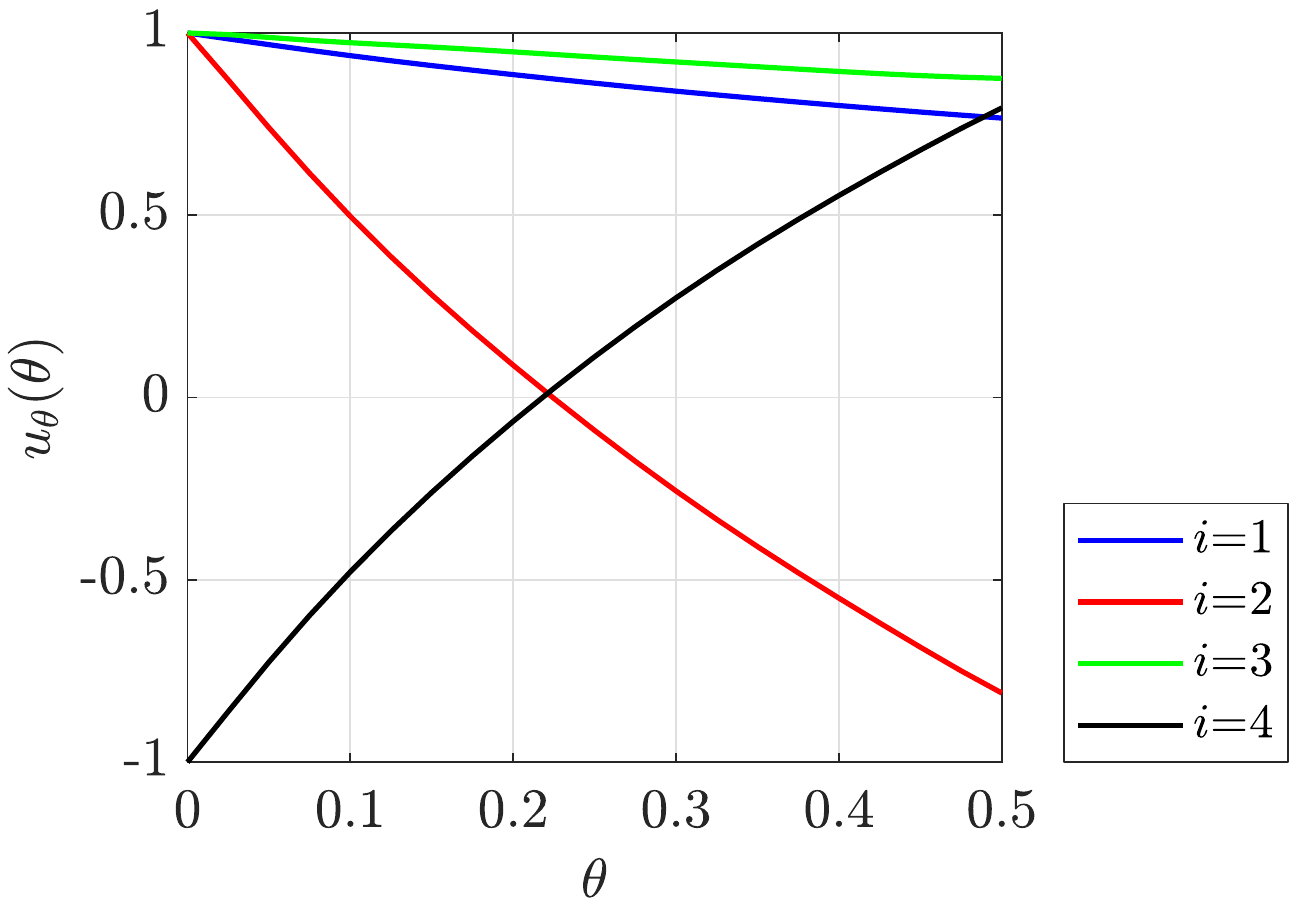}
	\end{subfigure}
	\caption{First four material $u_{\mu}(\mu)$ and geometric $u_{\theta}(\theta)$ modes of the generalized solution $\fU^{\texttt{PGD}}(\mu, \theta)$.}
	\label{fig:parametric_modes}
\end{figure}

\begin{figure}[tb!] 
	\small 
	\centering
	\includegraphics[width=1.\textwidth]{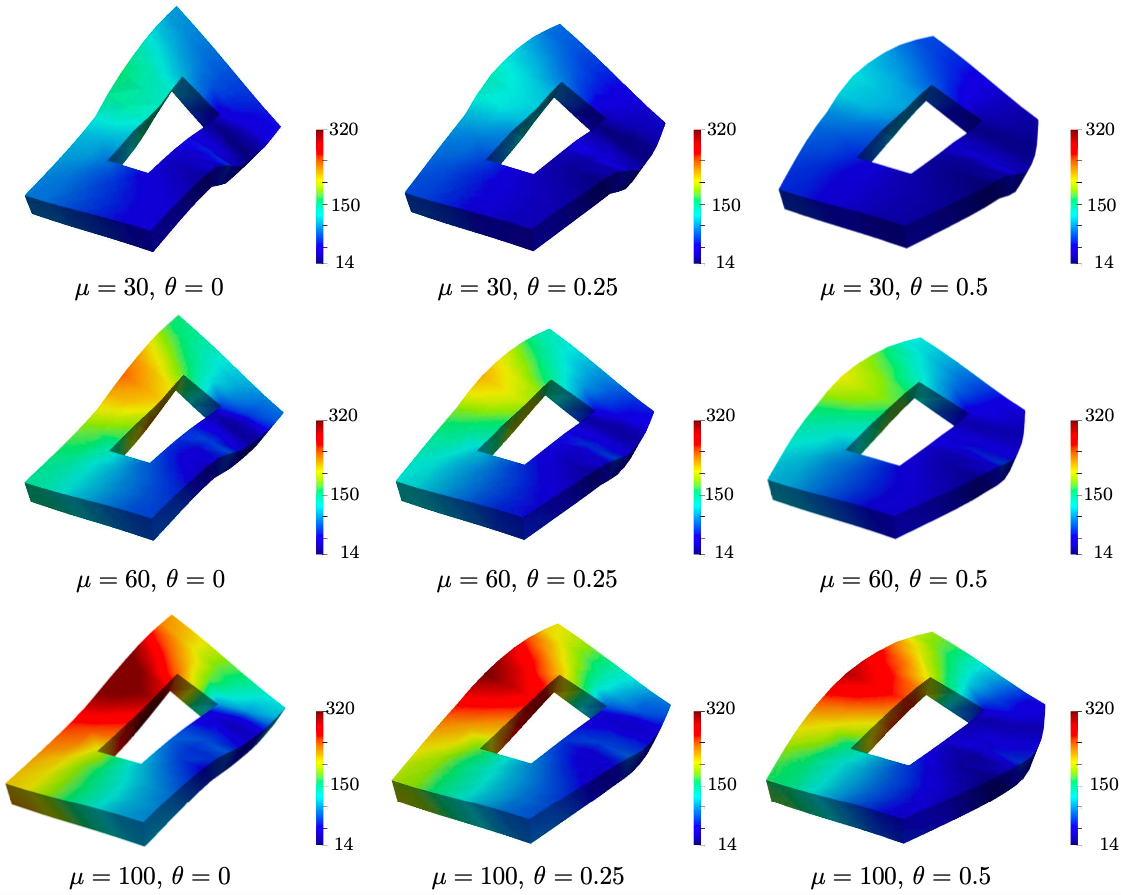}
	\caption{Particular cases of the generalized solution, showing the von Mises stress field, for nine choices of the parameters. The solutions are obtained in real-time after the PGD-IR is applied to compute the spatial and parametric modes.}
	\label{fig:donut_online}
\end{figure}

The first four normalized spatial modes are shown in Fig.~\ref{fig:spatialModes}, whereas the first four parametric modes are displayed in Fig.~\ref{fig:parametric_modes}. The spatial modes provide an illustration of the deformation induced by the four most relevant modes of the generalized solution. The parametric modes corresponding to the material illustrate that the four modes have the maximum contribution to the generalized solution for $\mu=10$. As the material property approaches the maximum value of $\mu=410$, the third and fourth mode have less influence on the solution. Finally, the modes corresponding to the geometric parameter have a more global character, proving the extra difficulty in solving geometrically parametrized problems.
In order to get a particularized solution for a chosen set of the parameters $(\bar \mu, \bar \theta)$, the correspondent function values $u^i_{\mu}(\bar \mu)$ and $u^i_{\theta}(\bar \theta)$ are evaluated for each PGD-mode $i$ and then multiplied by the correspondent spatial mode and amplitude.
Fig.~\ref{fig:donut_online} shows the particularized solutions in terms of deformed configuration and equivalent von Mises stress field for nine specific sets of parameters. The dominant character of the first spatial mode of Fig.~\ref{fig:spatialModes} can be clearly observed, whereas the magnitude of the stress highly depends on the parametric choice. Please remember that these particularized solutions were obtained in real-time during an online post-process step.

In order to validate the PGD results, the accuracy with respect to the full-order FE computations is measured as the relative error between the PGD and FE solutions in the $\mathcal{L}_2 (\Omega \times \mM_{\mu} \times \mM_{\theta})$ norm, that is 
\begin{equation}\label{eqRelErr}
\varepsilon_{\text{PGD}} = \left( \frac{ \int_{\mM_{\theta}} \int_{\mM_{\mu}} \int_{\Omega} \left(\bU^{\texttt{PGD}} - \bU^{\text{FE}} \right) \cdot \left(\bU^{\texttt{PGD}} - \bU^{\texttt{FE}} \right) d\Omega \: d\mu \: d\theta}{ \int_{\mM_{\theta}} \int_{\mM_{\mu}} \int_{\Omega}  \bU^{\texttt{FE}} \cdot \bU^{\texttt{FE}} \: d\Omega \: d\mu \: d\theta }	\right)^{1/2}.
\end{equation}
It is worth noting that to compute this error measure, the problem is solved by means of the standard FE method for each possible combination of the parameters, that is $n_{\mu} \times n_{\theta} = 21 \times 41 = 861$ FE simulations. 

\begin{figure}[tb!]
	\small
	\centering
	\includegraphics[width=.50\textwidth]{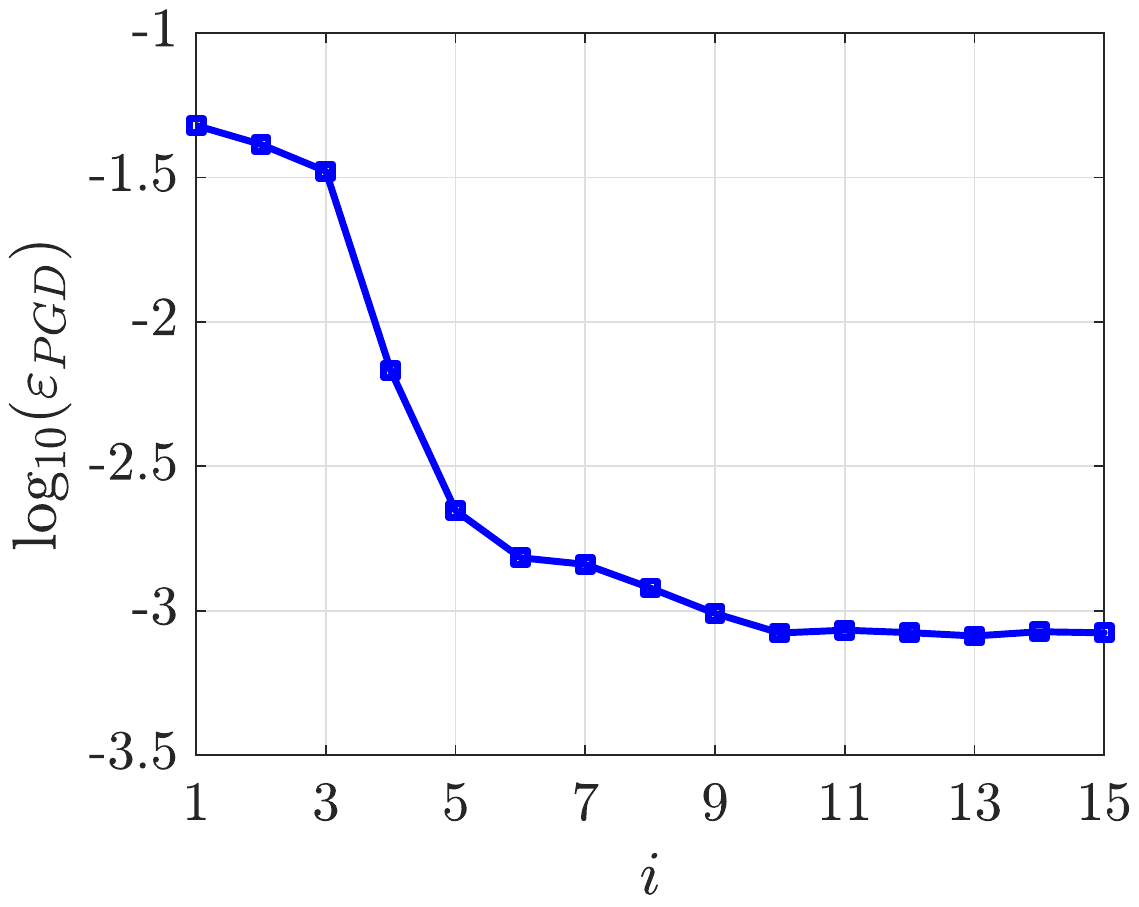}
	\caption{$\mathcal{L}_2 (\Omega \times \mM_{\mu} \times \mM_{\theta})$ norm of the difference between the PGD solution and the FE solution as a function of the number of PGD modes, $i$.}
	\label{fig:L2err}
\end{figure}

Fig.~\ref{fig:L2err} shows the evolution of the relative error with respect to the number of PGD modes. As expected, the level of accuracy increases as the number of modes increases, up to a user-defined tolerance, which in this case was chosen equal to $10^{-3}$. Note that the PGD solution converges to the desired tolerance with only nine PGD modes. 
An interesting advantage of the PGD method with respect to the standard FE method concerns the storage memory. In fact, the obtained PGD \textit{computational vademecum} needs \textasciitilde 74 KB of storage memory versus the \textasciitilde 6650 KB needed to store all the 861 full-order FE solutions. Computational time is not particularly significant in the PGD context. In fact, the main goal is to provide a method which is able to explore an arbitrary large parametric space with only one offline computation. Nevertheless, an interesting comparison is shown in Table~\ref{Table_iterations} where the number of iterations needed by the alternating direction scheme for the computations of each PGD mode is provided. As the cost of each iteration corresponds to the cost of a full-order FE simulation, the results in Table~\ref{Table_iterations} show that the cost of the proposed PGD-IR is equivalent to 161 full-order solutions, compared to the 861 full-order computations required by the standard FE approach.
\renewcommand{\arraystretch}{1.3}
\begin{table}[tb!]
  \centering
  \small
  \setlength\tabcolsep{1pt}
	\begin{tabular}{rrlllllrcc}
		\hline
		\multicolumn{1}{|l|}{\textbf{PGD mode}}         & \multicolumn{1}{r|}{\textbf{Mode 1}} &
		 \multicolumn{1}{c|}{\textbf{Mode 2}} &
		  \multicolumn{1}{c|}{\textbf{Mode 3}} &
		   \multicolumn{1}{c|}{\textbf{Mode 4}} &
		    \multicolumn{1}{c|}{\textbf{Mode 5}} &
		     \multicolumn{1}{c|}{\textbf{Mode 6}} &
		      \multicolumn{1}{c|}{\textbf{Mode 7}} &
		       \multicolumn{1}{c|}{\textbf{Mode 8}} &
		        \multicolumn{1}{c|}{\textbf{Mode 9}} \\ \hline
		\multicolumn{1}{|r|}{\textbf{N. of iterations}} & \multicolumn{1}{r|}{12}              & \multicolumn{1}{r|}{27}              & \multicolumn{1}{r|}{17}              & \multicolumn{1}{r|}{14}              & \multicolumn{1}{r|}{14}              & \multicolumn{1}{r|}{18}              & \multicolumn{1}{r|}{24}              & \multicolumn{1}{r|}{17}              & \multicolumn{1}{r|}{18}              \\ \hline
		& \textbf{}                            &                                      &                                      &                                      &                                      &                                      & \multicolumn{3}{r}{\textbf{Total n. of iterations = 161}}                                                                   
	\end{tabular}
	\caption{Total number of iterations performed by the alternating direction scheme to compute each PGD mode.}
	\label{Table_iterations}
\end{table}

\begin{figure}[tb!]
	\small 
	\centering
	\includegraphics[width=.55\textwidth]{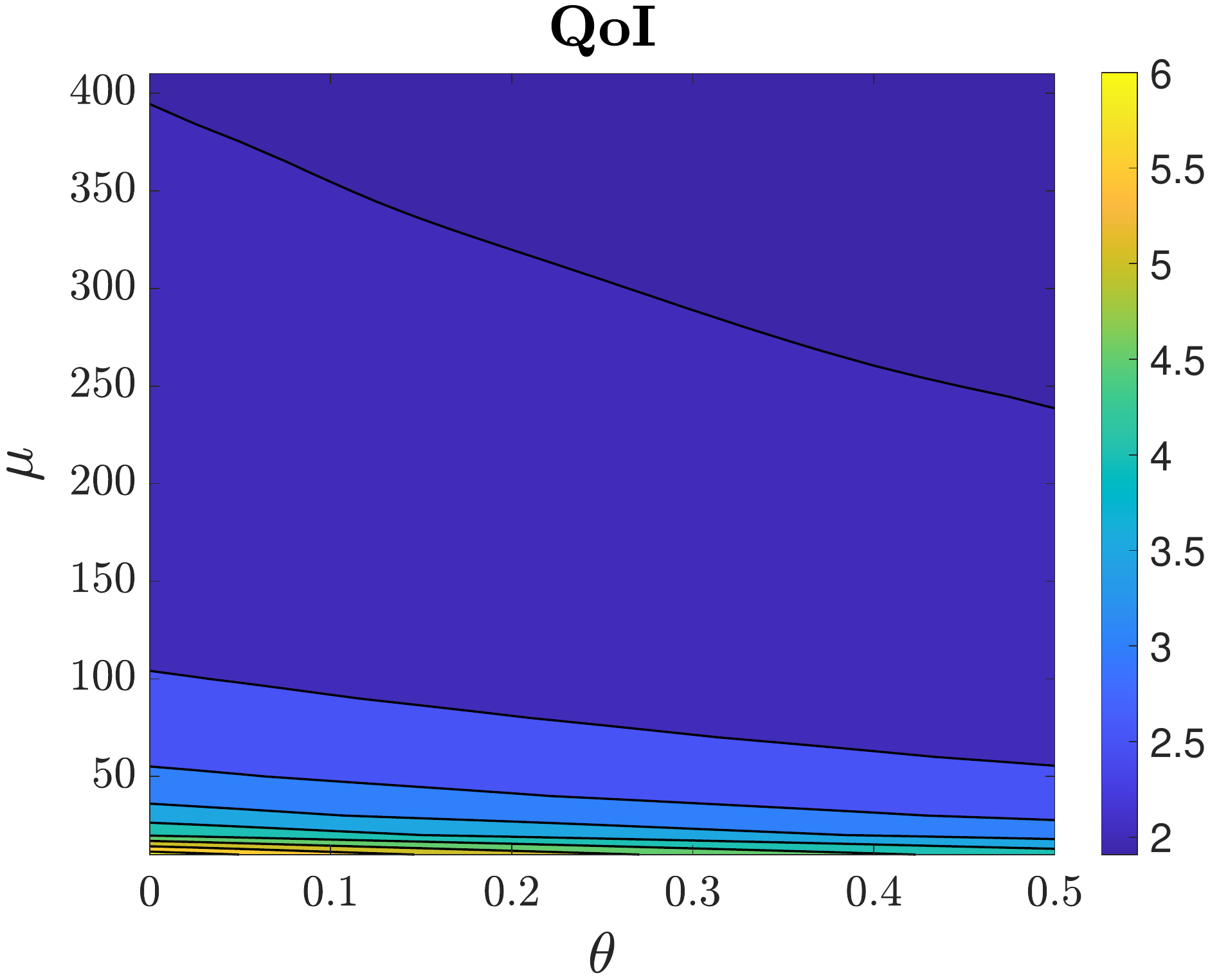}
	\caption{Variation of the QoI $\Delta U_{PQ}(z)$ with respect to the parameters $\mu$ and $\theta$.}
	\label{fig:QoI}
\end{figure}

Finally, a major advantage of computing a PGD \textit{computational vademecum} is the possibility to explore the design space and check, in real time, the effects of the design parameters on a predefined quantity of interest (QoI). As an example, the relative displacement in the $z$ direction, $\Delta \text{U}_{\text{PQ}} (z)$, of the points $P$ and $Q$ (see Fig.~\ref{fig:Donut}) is selected a QoI. The variation of the chosen QoI in the parametric space is depicted in Fig.~\ref{fig:QoI}. 

\subsection{Industrial application: dummy car test}\label{secCar}

\begin{figure}[tb!]
	\begin{subfigure}{.5\linewidth}
		\centering
		\includegraphics[width=\linewidth]{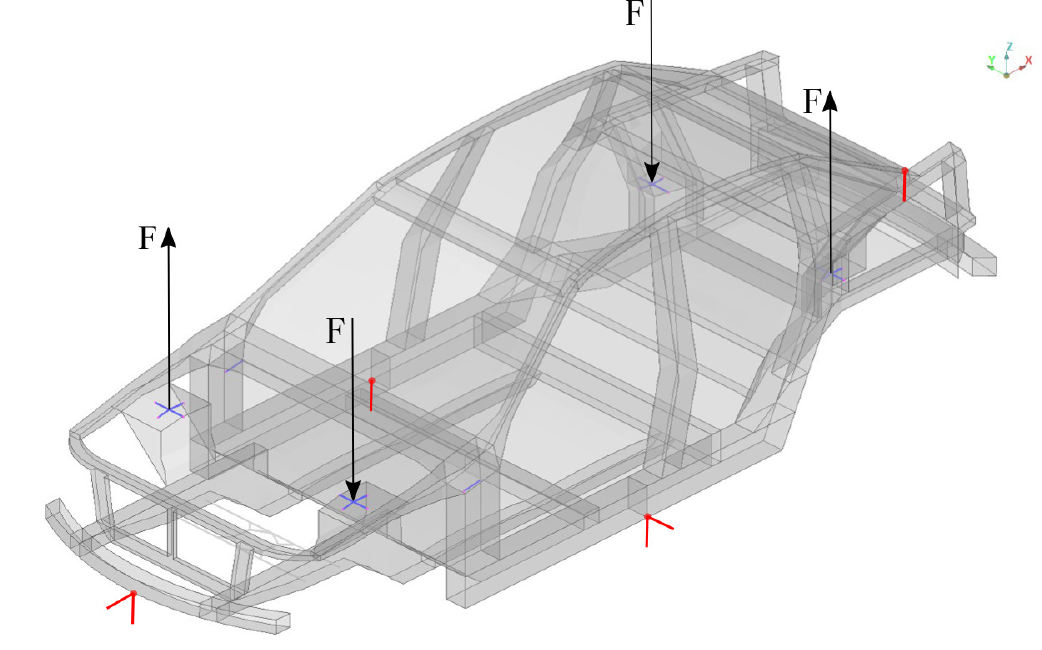}
	\end{subfigure}%
	\begin{subfigure}{.5\linewidth}
		\centering
		\includegraphics[width=\linewidth]{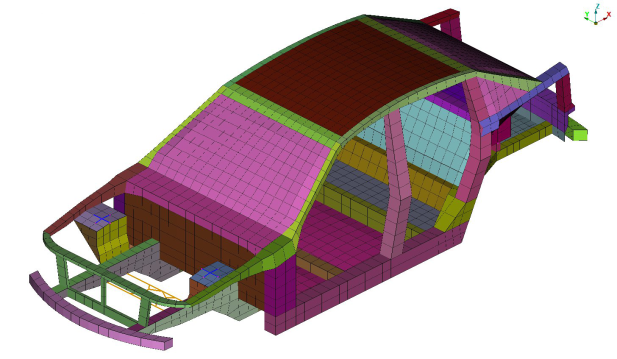}
	\end{subfigure}
	\caption{Geometry, load conditions (left) and mesh properties (right) of the BIW structure used for the static global torsional stiffness analysis.}
	\label{fig:dummycar}
\end{figure}

\begin{figure}[tb!]
	\centering
	\includegraphics[width=.55\textwidth]{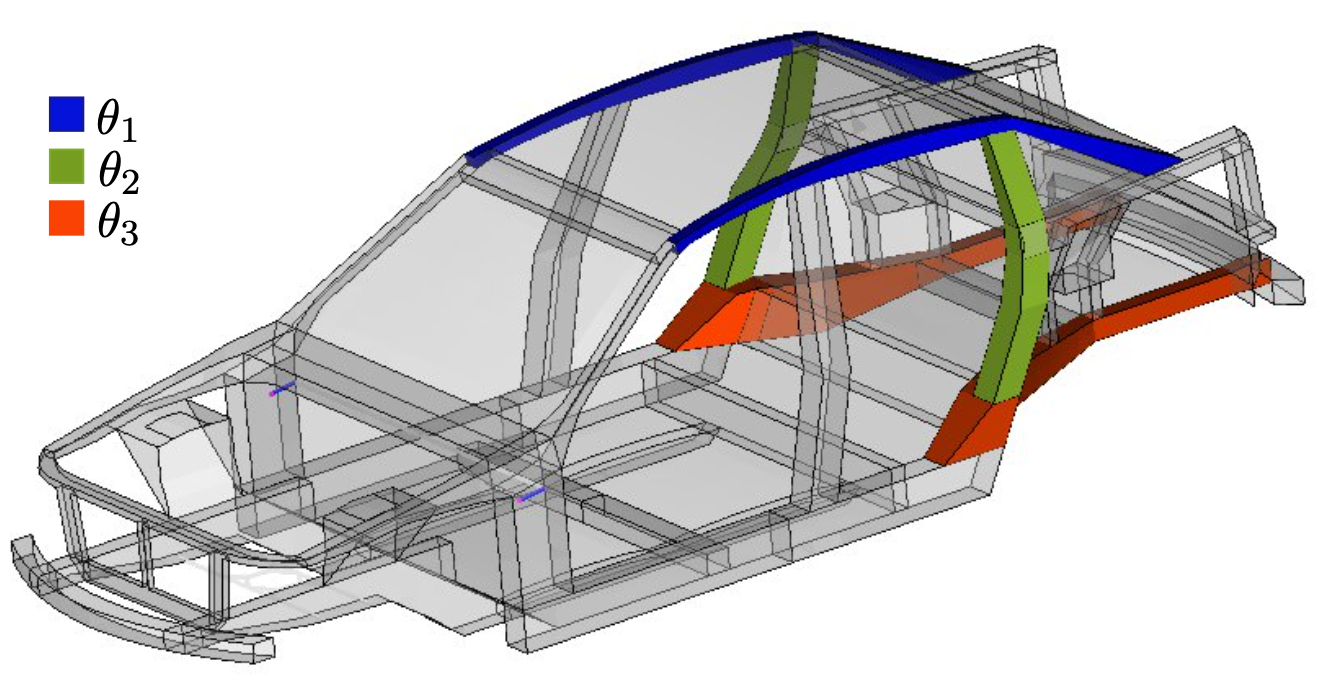}
	\caption{Geometry of the BIW with three car components highlighted. The parameters correspond to the thickness of each one of the components.}
	\label{fig:dummycar_param}
\end{figure}

\begin{figure}[tb!]
	\small 
	\centering
	\includegraphics[width=.55\textwidth]{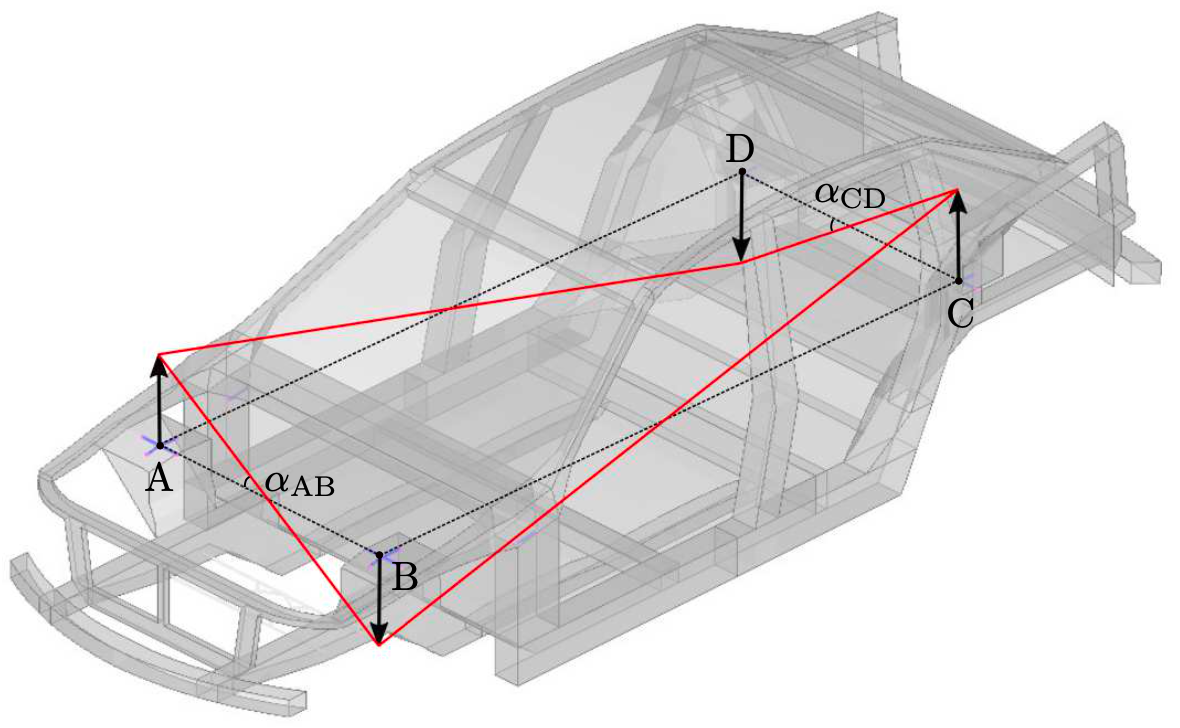}
	\caption{Illustration of the angles used to compute the ETS in Eq.~\eqref{eqETS}.}
	\label{fig:dummyCar_ETSdef}
\end{figure}

The PGD-IR method is now employed to solve a more realistic problem, which is the static global torsional stiffness analysis of the BIW structure of a generic car.
The geometry of the BIW is shown in Fig.~\ref{fig:dummycar}. Two couples of parallel and opposite forces are applied at the front and rear shock towers, such that two opposite torsional moments of magnitude 1 Nm are generated. The FE model, is discretized with isoparametric quadrilateral shell elements. The material is linear elastic and it is characterized by a Young's modulus $E = 207$ GPa, Poisson's ratio $\nu = 0.29$ and density $\rho = 7.82 \; kg/m^3$. 
In this example, the thickness of three car components highlighted in Fig.~\ref{fig:dummycar_param}, that usually play a role in the characterization of the global stiffness of the car, are introduced as extra coordinates of the problem. The three parameters are denoted by $\bmu = [\mu_1, \mu_2, \mu_3]^T$ and they vary in the intervals $\mM_j = [0.7, 1.5]$ mm, for $j = 1, 2, 3$. The three parametric domains are discretized with $n_1 = n_2 = n_3 = 9$ equidistant nodes.
The goal of this example is to demonstrate potential of the proposed PGD-IR approach, able to produce a generalized solution that enables a designed to check how the overall static stiffness of the vehicle is affected by any change of the introduced parameters. This is done by computing the equivalent torsional stiffness (ETS), which is defined as a function of the front and back twisting rotations of the car body when a torsion load is applied (see Fig.~\ref{fig:dummyCar_ETSdef}), namely
\begin{equation}\label{eqETS}
\text{ETS} = \frac{1}{\alpha_{\text{AB}}+ \alpha_{\text{CD}} } \times \frac{\pi}{180},
\end{equation}
where the two angles $\alpha_{\text{AB}}$ and $\alpha_{\text{CD}}$ are defined as
\begin{equation}\label{eqAlphaAB_CD}
\alpha_{\text{AB}} = \frac{\abs{U_z(A)} + \abs{U_z(B)}}{\norm{L_{\text{AB}}}}, \qquad 
\alpha_{\text{CD}} = \frac{\abs{U_z(C)} + \abs{U_z(D)}}{\norm{L_{\text{CD}}}}. 
\end{equation}
In the above expressions, $U_z(P)$ denotes the displacement in the $z$ direction at point P and $L_{\text{PQ}}$ denotes the distance between the points P and Q.

The proposed PGD-IR approach is employed following the same procedure described in the previous example. In the preprocess stage, the commercial FE package MSC-Nastran is now used to sample the parametric input matrices. A script was prepared to automatically produce a new Nastran input file (.bdf and .dat files) for each possible combination of the parameters. The generated files were then read by the Nastran solver, where the matrices were assembled (without solving the problem) and stored in a plain text format. Afterwards, the matrices were read by the developed Matlab routine to be expressed in the required separated form, namely
\begin{equation}\label{eqK_dummycar}
\fK^{\texttt{PGD}}(\bmu) = \sum_{i=1}^{n_{\text{tot}}} \fK^i \: \prod_{j = 1}^{3} k_j^i(\mu_j),
\end{equation}
with $n_{\text{tot}} = n_1 \times n_2 \times n_3 = 729$, while $\fK^i$ and $k_j^i(\mu_j)$ being defined as described in Section~\ref{secAlgSep}. The mass matrix was obtained in a similar fashion and the separated force vector was computed. Finally, the separated expression of both the parametric stiffness and mass is compressed to minimize the number of terms in the separated form.

As already shown in the previous example, once the input data $\fK^{\texttt{PGD}}(\bmu), \fM^{\texttt{PGD}}(\bmu)$ and $\fF^{\texttt{PGD}}(\bmu)$
are pre-processed and expressed in a separated form, the cascade scheme of the encapsulated PGD is used to sequentially solve the three steps involved in the developed PGD-IR approach, such that the final solution is obtained as
\begin{equation}\label{eqU_dummycar}
\fU(\bmu) ^{\texttt{PGD}}= \sum_{i=1}^{N_{\bU}} \beta_{\bU} \fU^i \: \prod_{j = 1}^{3} u_j^i(\mu_j).
\end{equation}

\begin{figure}[tb!]
	\small
	\centering
	\includegraphics[width=.5\textwidth]{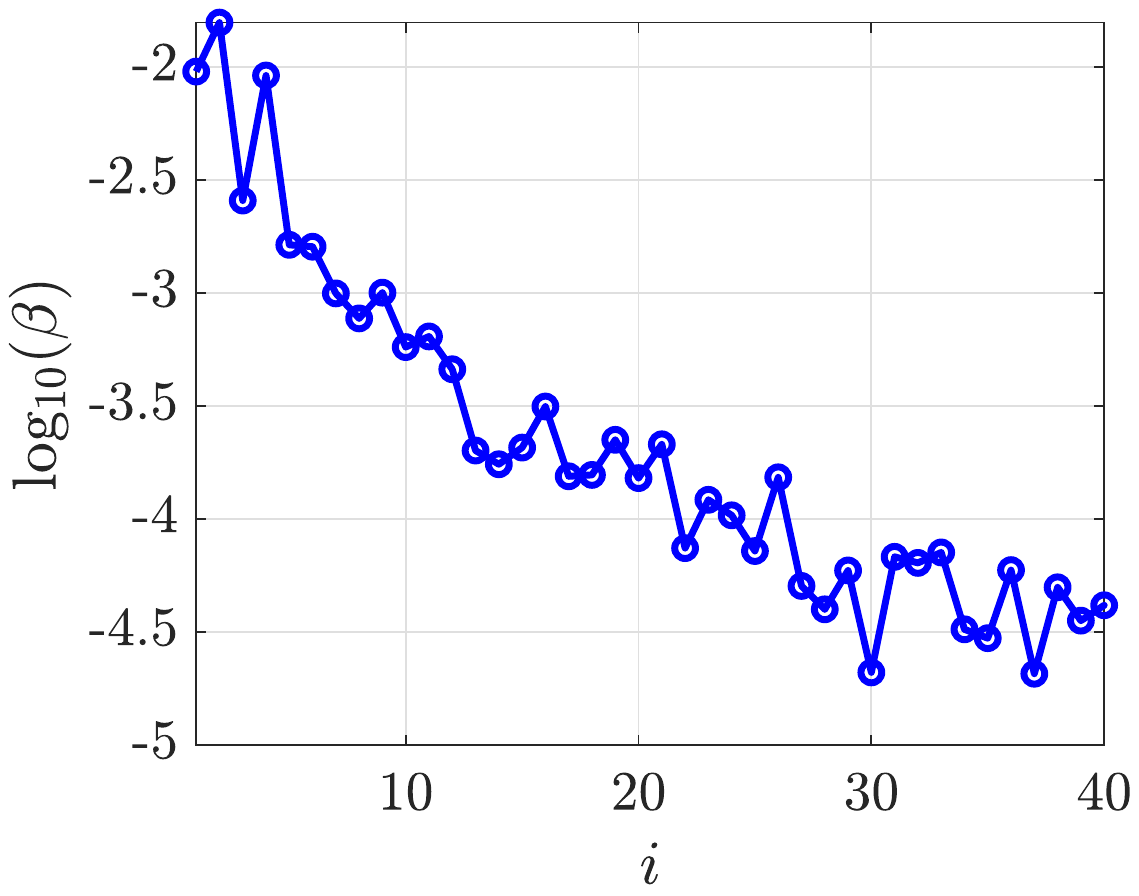}
	\caption{Evolution of the amplitude of the PGD modes, $\beta^i$, of the solution $\fU^{\texttt{PGD}}(\bmu)$}
	\label{fig:dummycar_amplitudes}
\end{figure}

The amplitude of the PGD modes is shown in Fig.~\ref{fig:dummycar_amplitudes}. In this example, with a more complex geometry and a larger number of parameters, it can be observed that more modes are required to produce an accurate description of the multi-dimensional solution. With 21 modes the amplitude of the modes is approximately two orders of magnitude lower than the amplitude of the first mode. 

\begin{figure}[tb!]
	\begin{subfigure}{.24\linewidth}
		\centering
		\includegraphics[width=1\linewidth]{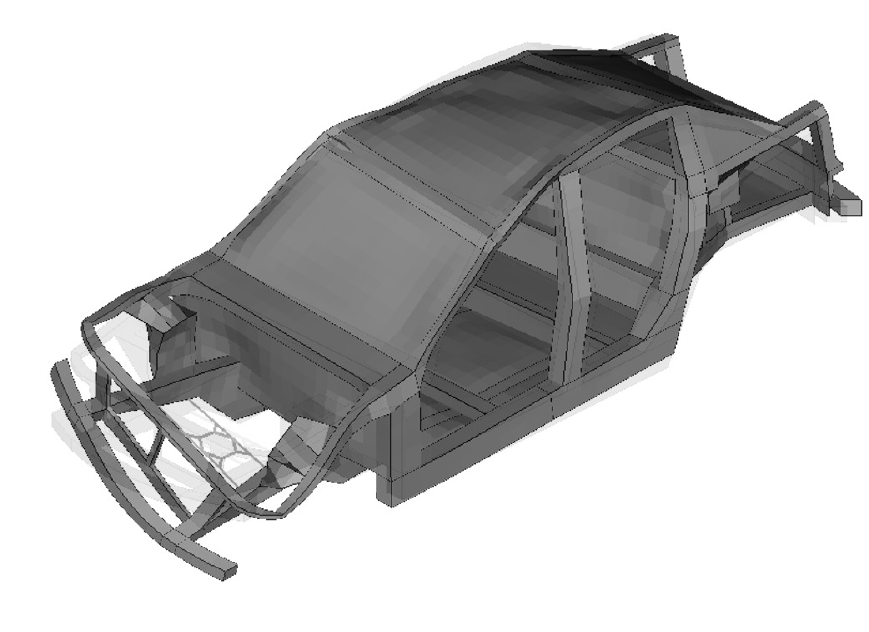}
		\caption{Mode 1}
		\label{fig:spatial1}
	\end{subfigure}%
	\begin{subfigure}{.24\linewidth}
		\centering
		\includegraphics[width=1\linewidth]{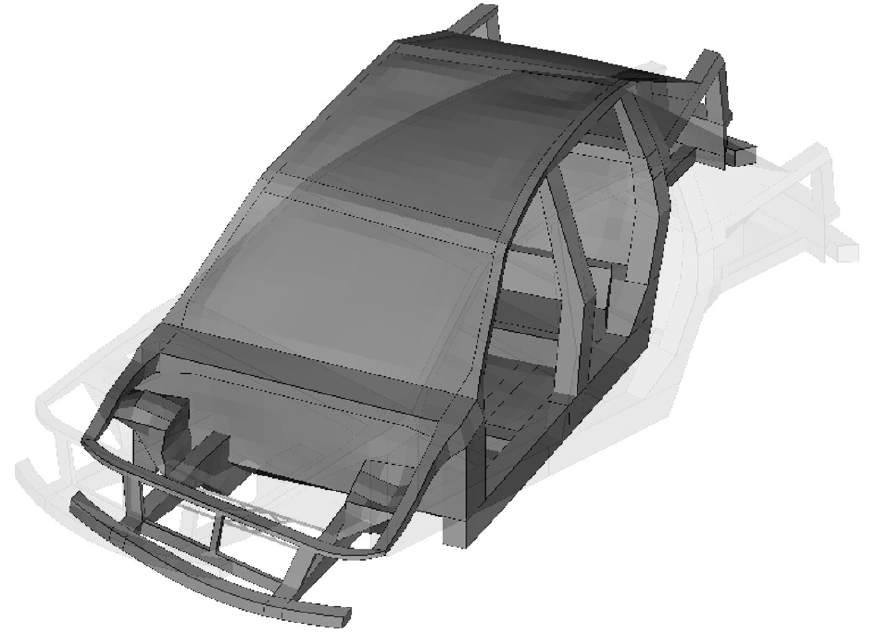}
		\caption{Mode 2}
		\label{fig:spatial2}
	\end{subfigure}
	\begin{subfigure}{.24\linewidth}
		\centering
		\includegraphics[width=1\linewidth]{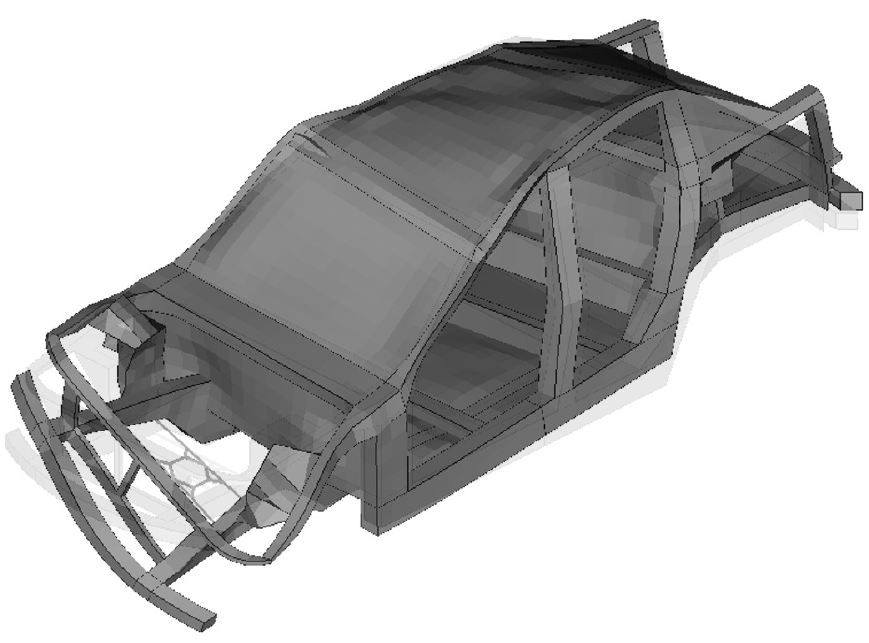}
		\caption{Mode 3}
		\label{fig:spatial3}
	\end{subfigure}%
	\begin{subfigure}{.24\linewidth}
		\centering
		\includegraphics[width=1\linewidth]{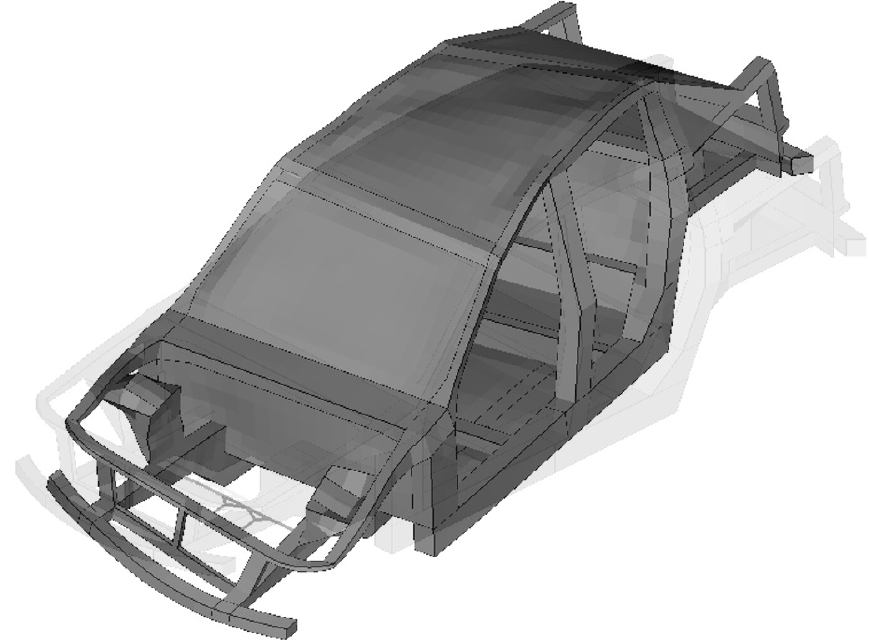}
		\caption{Mode 4}
		\label{fig:spatial4}
	\end{subfigure}
	\caption{ First four spatial modes of the solution generalized solution $\fU^{\texttt{PGD}}(\bmu)$.}
	\label{fig:dummycar_spatialModes}
\end{figure}

\begin{figure}[tb!]
	\begin{subfigure}{.45\linewidth}
		\centering
		\includegraphics[width=1\linewidth]{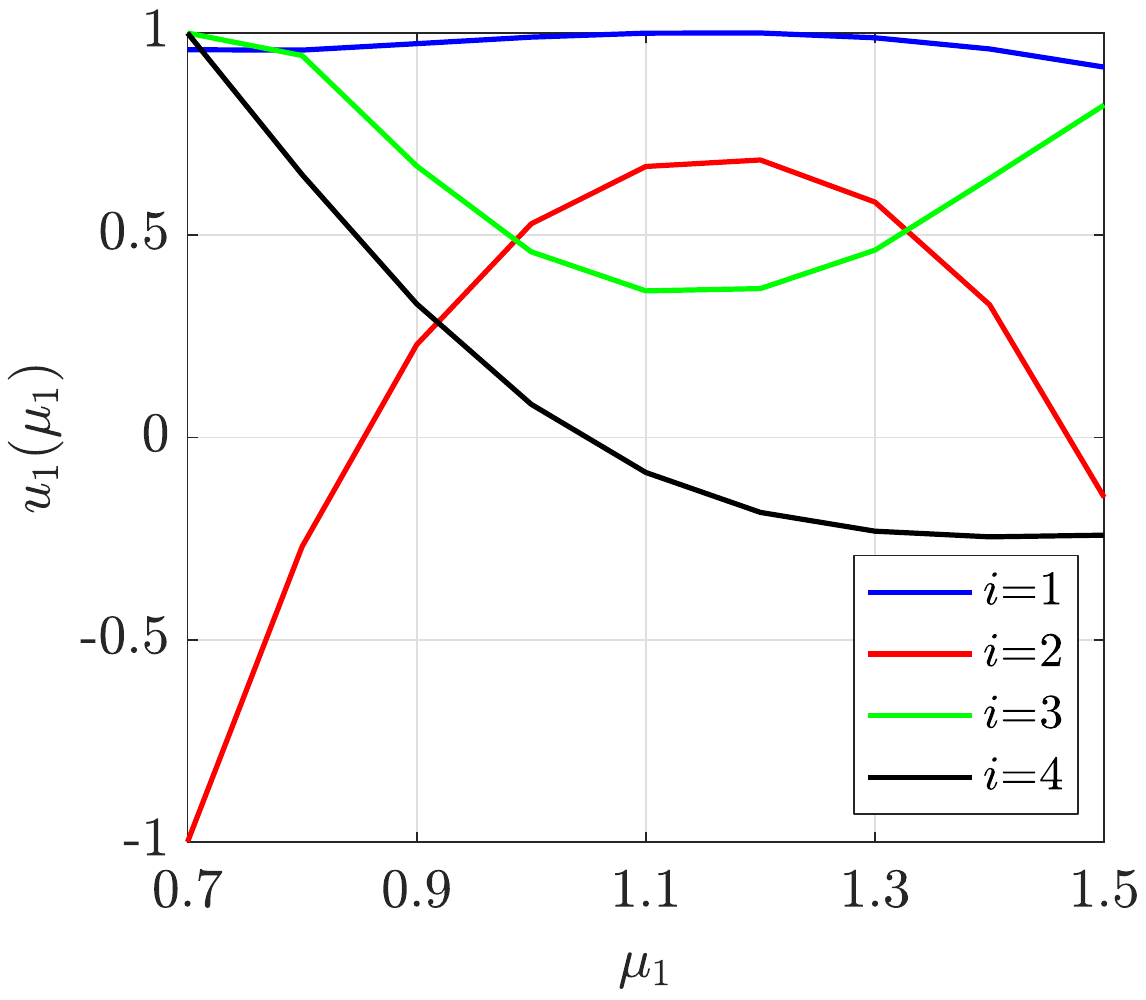}
		\caption{$u_1(\mu_1)$}
		\label{fig:f_mu}
	\end{subfigure}%
	\begin{subfigure}{.45\linewidth}
		\centering
		\includegraphics[width=1\linewidth]{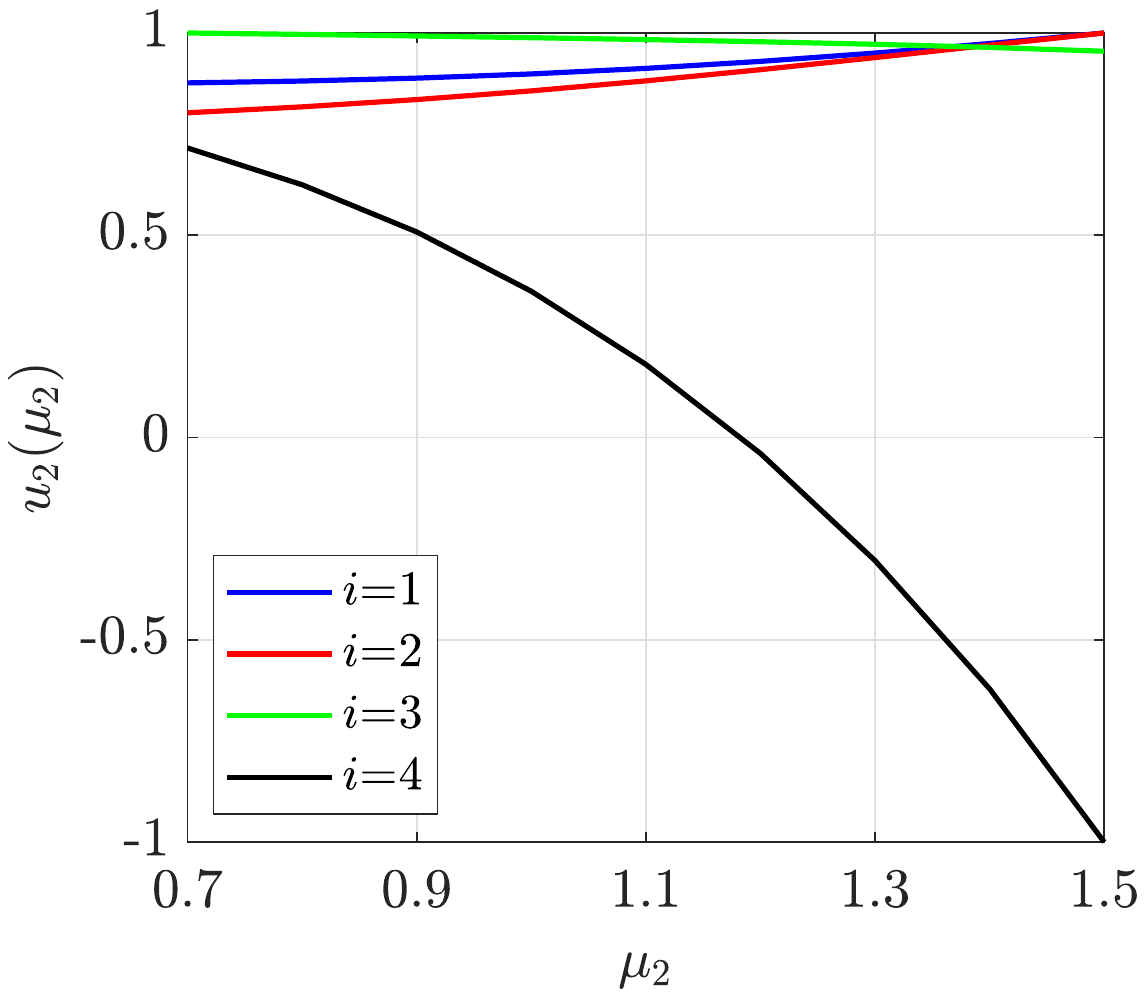}
		\caption{$u_2(\mu_2)$}
		\label{fig:f_theta}
	\end{subfigure}
	\begin{subfigure}{.45\linewidth}
		\centering
		\includegraphics[width=1\linewidth]{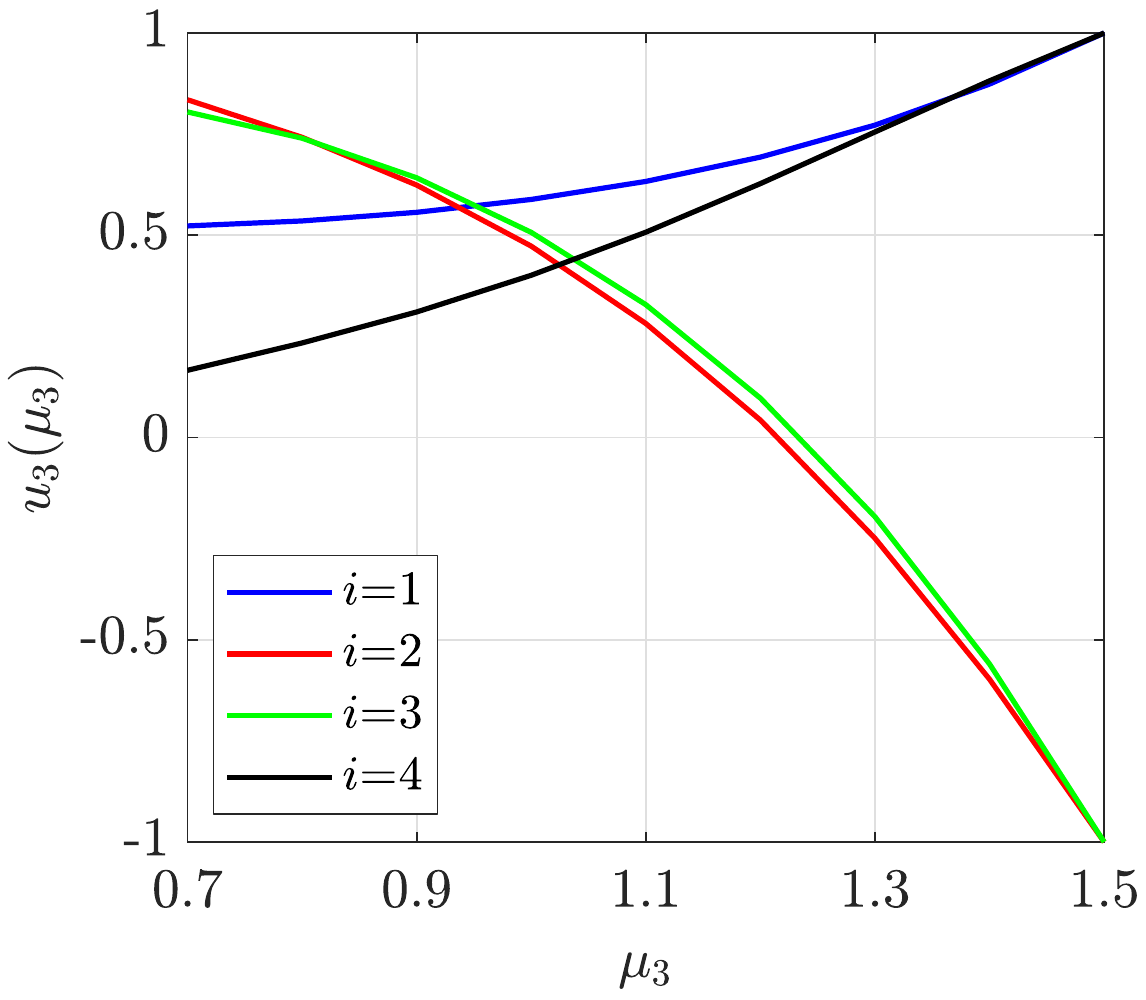}
		\caption{$u_3(\mu_3)$}
		\label{fig:f_theta}
	\end{subfigure}
	\caption{First four parametric modes of the generalized solution $\fU^{\texttt{PGD}}(\bmu)$.}
	\label{fig:dummycar_paramModes}
\end{figure}

The first four spatial modes are depicted in Fig.~\ref{fig:dummycar_spatialModes} amplified by a factor of  \textasciitilde1000 and Fig. ~\ref{fig:dummycar_paramModes} shows the first four normalized parametric functions.

To illustrate the full potential of the proposed PGD-IR approach, the application of the PGD-IR in a multi-objective optimization process is considered. The goal is to find the combination of parameters that maximize the ETS, while minimizing the mass of the three car components considered in this example. To this end, two objective functions are considered, namely
\begin{align}
\begin{cases}\label{eqObjFun}
g_1(\bmu) &= \rho \: (\mu_1 A_1 + \mu_2 A_2 + \mu_3 A_3), \\
g_2 (\bmu) &= \text{ETS}(\bmu),
\end{cases}
\end{align}
where $g_1(\bmu)$ represents the mass of the material needed to manufacture the three car components, equal to the product of the material density $\rho$ and the parametric volume. The latter is given by the sum of the products between the car components areas $(A_1, A_2, A_3)$ times their variable thicknesses $(\mu_1, \mu_2, \mu_3)$. Clearly, this quantity is strictly related to the production cost. The objective function $g_2(\bmu)$ represents the parametric ETS defined in Eq.~\eqref{eqETS}.

With the computed generalized solution $\fU^{\texttt{PGD}}(\bmu)$, an evaluation of the objective functions within a multi-objective optimization process only requires the particularization of the solution for a given set of the parameters. With the proposed PGD-IR approach, this evaluation can be performed in real time, making the overall cost of the optimization stage almost negligible. This is in contrast with a traditional approach where each evaluation of the objective functions require the assembly and solution of a new FE system of equations. 

In this example, the optimization problem was solved by means of the \texttt{gamultiobj} function available in the Global optimization Toolbox released by Matlab. The function is able to find the Pareto front of multiple objective functions using a genetic algorithm. A Pareto front is a set of optimal points in the parametric space that represent a trade-off between the objective functions. More specifically, a point is considered optimal if no objective can be improved without sacrificing at least one other objective. 

\begin{figure}[tb!]
	\begin{subfigure}{.5\linewidth}
		\centering
		\includegraphics[width=1.\linewidth]{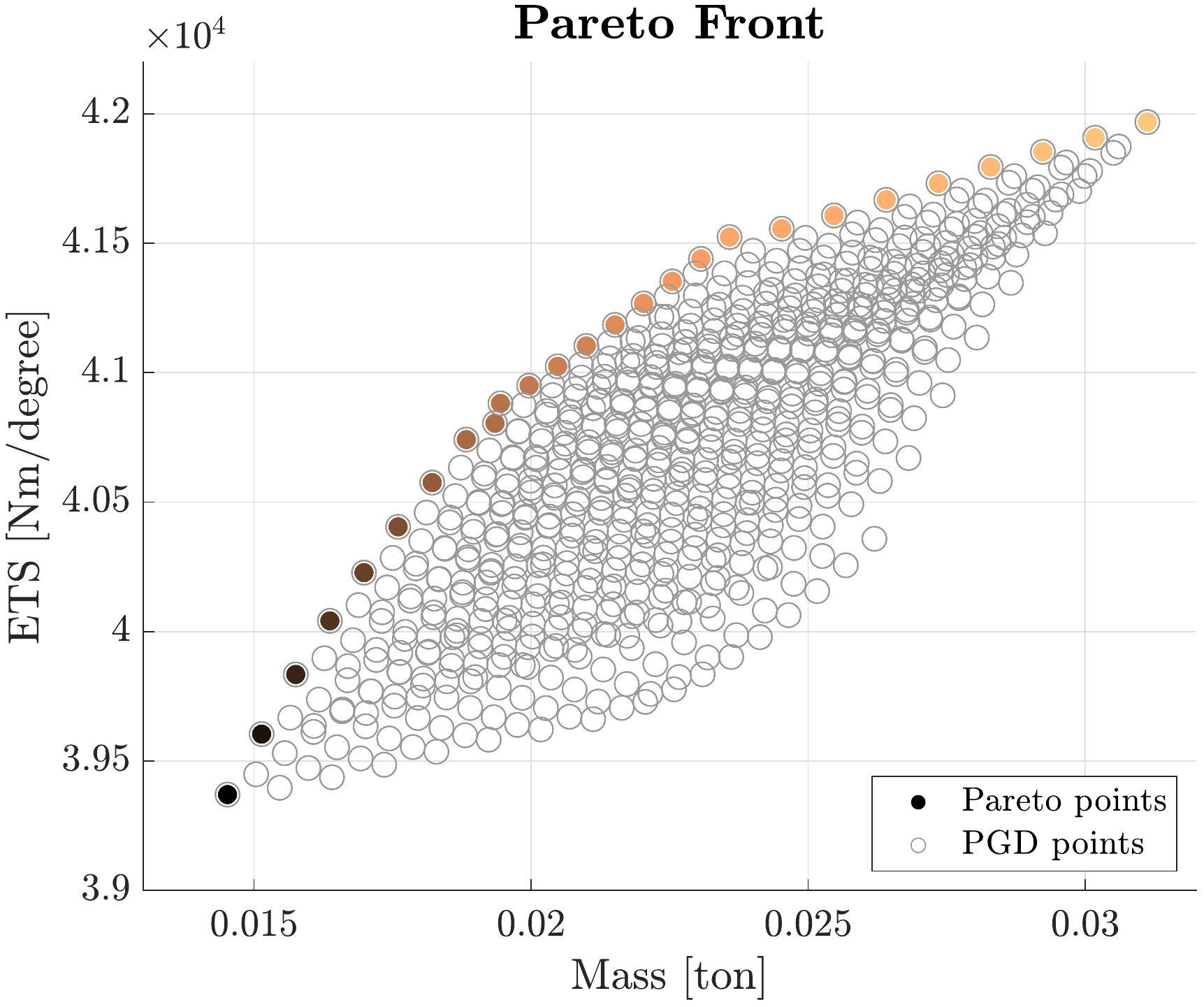}
		\label{fig:ParetoFront}
	\end{subfigure}
	\begin{subfigure}{.5\linewidth}
		\centering
		\includegraphics[width=1.\linewidth]{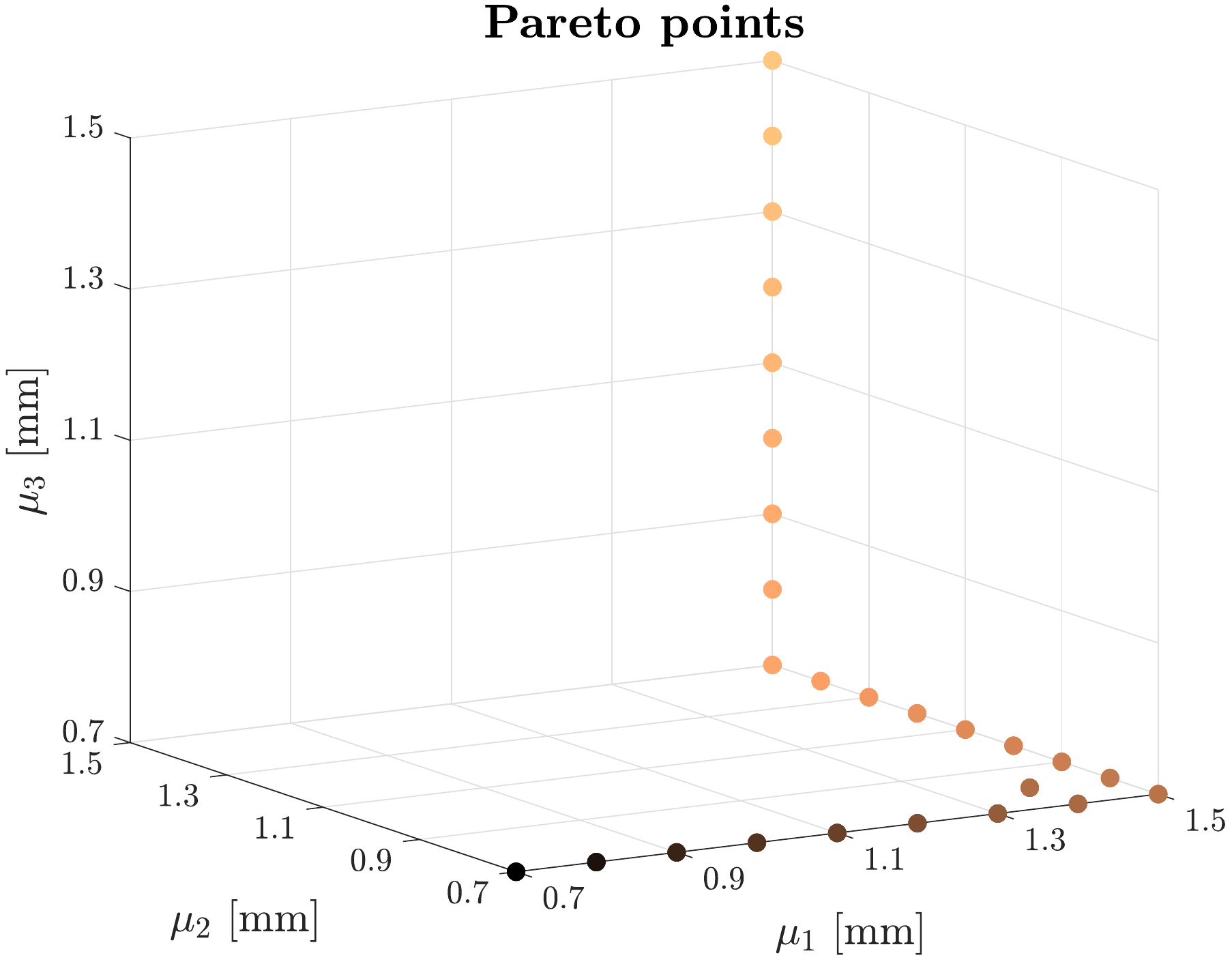}
		\label{fig:ParetoPoints}
	\end{subfigure}
	\caption{Multi-objective optimization showing the Pareto front as a function of the objectives (left) and the PGD parameters (right).}
	\label{fig:Pareto}
\end{figure}

Fig.~\ref{fig:Pareto} shows the Pareto front in terms of the objective functions as well as the whole range of configurations that result from the parametric space $\mM_1 \times \mM_2 \times \mM_3$. It is important to note that the optimization study allows to significantly reduce the range of solutions to be considered by a designer in the decision-making process. In fact, the variability in the two quantities of interest (mass and ETS), induced by the three parameters and calculated by PGD for all the possible combinations (PGD points in the left Fig.~\ref{fig:Pareto}), is much larger than the number of points belonging to the Pareto front. Fig.~\ref{fig:Pareto} (right) plots the Pareto points in the parametric space $(\mu_1, \mu_2, \mu_3)$, where the correspondence to the Pareto front is described by same colors. In this example, the Pareto front was computed by assigning the same weight to the objective functions. Nevertheless, it is straightforward to obtain other fronts if the user wants to put more emphasis on one of the objective functions.
\section{Conclusions}\label{secConcl}

A nonintrusive algebraic PGD approach combined with the IR method for the solution of unconstrained problems being characterized by material and geometric parameters has been presented. The developed solver makes use of the \textit{Encapsulated PGD Toolbox} developed by D\'iez et al.~\cite{diez2019encapsulated}, which enables to perform algebraic operations for multidimensional data and allows to solve sequentially the three parametric problems required by the IR method. 

An algebraic approach has been proposed to deal with geometric parameters by morphing a mesh generated for a reference configuration. The proposed method acts as a \textit{black-box}, such that a nonintrusive interaction with commercial FE packages is possible. 

Two numerical examples are used to underline the main properties of the method. The first example considers an academic test case with one material and one geometric parameter. The ability to compute a \textit{computational vademecum} is shown and the accuracy of the generalized solution is measured by comparing the PGD solution to a set of stadard FE full-order solutions. It it shown that the proposed PGD-IR approach is able save almost the 99\% of storage memory, requiring only the 20\% of computational time needed by the FE method to solve the problem for every possible set of parameters. The second problem involves an industrial application for the static global stiffness analysis of a BIW structure of a car characterized by three parameters. This example shows the potential of the proposed PGD-IR approach and its ability to be integrated with a commercial FE package, such as MSC-Nastran. Finally, a multi-objective optimization was performed in order to show how the proposed approach can represent an important support to designers during the decision-making process. With the developed technique it is to possible to produce a \textit{computational vademecum} displayed in a portable device to support the design engineers in the decision-making by evaluating in real time the impact of certain parameters on the global response of the structure. 

\section*{Acknowledgments}
This project is part of the Marie Skłodowska-Curie ITN-EJD ProTechTion funded by the European Union Horizon 2020 research and innovation program with grant number 764636. The work of Fabiola Cavaliere, Sergio Zlotnik and Pedro Díez is partially supported by the Spanish Ministry of Economy and Competitiveness, Spain (Grant number: DPI2017-85139-C2-2-R) and by the Generalitat de Catalunya, Spain (Grant number: 2017-SGR-1278). Ruben Sevilla also acknowledges the support of the Engineering and Physical Sciences Research Council (Grant number: EP/P033997/1).
\bibliographystyle{elsarticle-num}
\bibliography{wileyNJD_AMA}
\appendix

\section{Analytical approach to separate input quantities}\label{app1}

The analytical technique follows the standard isoparametric concept widely used in FE formulations. This implies that a mapping function ${\Psi}(\theta)$ which transforms the reference domain $\Omega$ into the geometrically parametrised domain ${\Omega}(\theta)$ has to be defined, such that
\begin{align}
\Psi(\theta) :  \; &\Omega \rightarrow \Omega(\theta) \\
& \mathbf{X} \mapsto  \mathbf{x}  = \Psi(\mathbf{X}, \theta),
\end{align}
where $\mathbf{X}$ represents the coordinate system associated to the reference domain $\Omega$, while $\mathbf{x}$ describes the modified domain $\Omega(\theta)$. According to the standard procedure, in order to transform the integrals involved in the weak formulation from the parametrized domain to the reference one, the Jacobian matrix $\mathbf{J}_{\Psi}(\theta) = \partial{\mathbf{x}}/\partial{\mathbf{X}}$ associated to the mapping $\Psi(\theta)$ has to be introduced. Then, the discretized definition of the stiffness and mass matrices at the reference element level $\Omega_e$ (already defined in Sec. \ref{secDiscr}) becomes
\begin{equation}\label{eqStiffElemParam}
\mathbf{K}^e = \int_{ \Omega_e} \T{\mathbf{B}} \T{{\mathbf{J}_{\Psi}}(\theta)^{-}} \mathbfcal{C} \: {\mathbf{J}_{\Psi}}(\theta)^{-1} \mathbf{B} \;  \text{det}(\mathbf{J}_{\Psi}(\theta)) \; \text{d}\Omega,
\end{equation}
\begin{equation}\label{eqMassElemParam}
\mathbf{M}^e = \int_{ \Omega_e} \T{\mathbf{N}}  \mathbf{N} \;  \text{det}(\mathbf{J}_{\Psi}(\theta)) \; \text{d}\Omega.
\end{equation}

The modification of the stiffness matrix formulation caused by the introduction of the inverse of the Jacobian matrix $\mathbf{J}_{\Psi}^{-1}$, leads to the first important limitation of the analytical method. In fact, it is well known that even if the mapping function and correspondent Jacobian can be written in a separated form, the inverse matrix $\mathbf{J}_{\Psi}^{-1}(\theta)$ is in general not separable due to the presence of $\det(\mathbf{J}_{\Psi}(\theta))$ at the denominator. As a consequence, an explicit dependency of the stiffness matrix on the geometric parameter cannot be found and other methods should be employed to find a separated expression of it. This difficulty is discussed in detail in~\cite{sevilla2020solution}. An alternative mixed formulation has also been recently considered to circumvent this difficulty in a discontinuous Galerkin framework~\cite{sevilla2020hybridisable}. 
\end{document}